\documentclass[aip, reprint]{revtex4-1}
\usepackage{graphicx}
\usepackage{color}
\usepackage{amssymb,amsmath,amsfonts,mathtools,amsthm}
\usepackage{array,multirow,hhline}
\usepackage{makecell}
\usepackage{empheq}
\usepackage{stmaryrd}

\usepackage{hyperref}
\usepackage{hyperref}
\hypersetup{
    colorlinks,
    citecolor=blue,
    filecolor=blue,
    linkcolor=blue,
    urlcolor=blue
}

\newtheorem*{conj}{Conjecture}

\begin{document}

\title{Topological Entropy of Surface Braids and Maximally Efficient Mixing}
%Might want to position the idea of surface braids a bit more prominently in the paper:
%%"Surface Braids and Mixing Efficiency"
%%Graph Generated Surface Braids and Maximal Mixing Efficiency
%%Topological Entropy of Graph-Generated Surface Braids and Maximally Efficient Mixing
%%Original Title:  Braids on Lattice Graphs and Maximally Efficient Mixing

\author{Spencer Ambrose Smith}
\author{Sierra Dunn}
\affiliation{Mount Holyoke College, South Hadley, MA 01075}

\date{\today}
\begin{abstract}

The deep connections between braids and dynamics by way of the Nielsen-Thurston classification theorem have led to a wide range of practical applications.  Braids have been used to detect coherent structures and mixing regions in oceanic flows, drive the design of industrial mixing machines, contextualize the evolution of taffy pullers, and characterize the chaotic motion of topological defects in active nematics.  Mixing plays a central role in each of these examples, and the braids naturally associated with each system come equipped with a useful measure of mixing efficiency, the topological entropy per operation (TEPO).  This motivates the following questions.  What is the maximum mixing efficiency for braids, and what braids realize this?  The answer depends on how we define braids.  For the standard Artin presentation, well-known braids with mixing efficiencies related to the golden and silver ratios have been proven to be maximal.  However, it is fruitful to consider surface braids, a natural generalization of braids, with presentations constructed from Artin-like braid generators on embedded graphs.  In this work, we introduce an efficient and elegant algorithm for finding the topological entropy and TEPO of surface braids on any pairing of orientable surface and planar embeddable graph.  Of the myriad possible graphs and surfaces, graphs that can be embedded in $\mathbb{R}^2$ as a lattice are a simple, highly symmetric choice, and the braids that result more naturally model the motion of points on the plane.  We extensively search for a maximum mixing efficiency braid on planar lattice graphs and examine a novel candidate braid, which we conjecture to have this maximal property.

%Include 

\end{abstract}

\maketitle
%\tableofcontents
\newpage

\section{Introduction}
\label{Sec:Intro}

Much of the allure of braid theory comes from the coexistence of deep, multi-disciplinary mathematical problems alongside numerous practical applications.  This is particularly evident with aspects of braid theory that can be viewed from the perspective of dynamical systems.  

Fluid dynamics has been a large source of inspiration for applications, starting with the work of Boyland, Aref, and Stremler\cite{MR1742169}, and continuing with the work of Thiffeault~\cite{MR2464304,MR2317905,BraidsEntPartTraj,PhysRevE.73.036311}.  In this context, the mixing of a 2D fluid (a fluid surface, or a bulk fluid where disparate characteristic time-scales result in effective 2D motion) is encoded in the braid formed from the space-time trajectories of stirring rods or passively advected particles.  For example, the trajectories of data-collection buoys can help reveal regions of the ocean in which the flow enhances mixing or help identify coherent structures in which mixing is minimized~\cite{ALLSHOUSE201295,PhysRevFluids.5.054504}.  These structures~\cite{LCShaller} play a central role in the transport of nutrients, oxygenation, temperature, and pollutants.  

More broadly, ideas from braid theory have been used to characterize the periodic orbits of point vortex motion~\cite{BOYLAND200369,smith2015point}, to investigate mixing in lid-driven cavity flow~\cite{doi:10.1063/1.2772881} and channel flow~ \cite{doi:10.1063/1.3076247}, and have been extended to include almost-invariant sets~ \cite{PhysRevLett.106.114101,doi:10.1063/1.4768666}.

In many industrial applications, mixing is a desirable outcome.  A whole class of machines, those which create mixing using the motion of stirring rods embedded in the fluid, are designed based off of braids which maximize mixing efficiency~\cite{MR2861264}.  These are especially useful for fluids in the Stokes regime (low Reynolds number, e.g. high viscosity), where turbulent mixing is negligible~\cite{ReynoldMixSmith}.

A more whimsical application looks at the evolution of taffy-pulling machines~\cite{TaffyThiffeault}.  Taffy is a soft candy made from repeatedly stretching melted sugar; in this process air is folded in,  resulting in a softer, more desirable texture.  Curiously, many of the machines designed in the late $19^{th}$ and early $20^{th}$ centuries have movements consistent with braids of high mixing efficiency.

As a final application, braids have been used to analyze the dynamics of topological defects in an active nematic microtubule system~\cite{ActiveNemTopChaos2019}.  Here, biologically-derived microtubule bundles confined to 2D define a director field; defects in this field are topologically protected and can be treated as particles.  As energy is injected into the system (at small scales, using ATP), the bundles extend in length causing a flow on the large scale.  Interestingly, the extensile dynamics ensure that the braids formed by the motion of topological defects have very particular mixing efficiencies.

In each of these examples the movement of a set of point-like objects, whether they are oceanic buoys, mixing rods, or topological defects, define trajectories which wind about one another to form mathematical braids.  In some applications the braids are passively formed, and we use them as a way of discretely encoding the salient aspects of the flow's complexity~\cite{MR3456018}.  In other applications the braids are actively created, and we use them to impose a lower bound on the flow's complexity.  In both cases, we appeal to one of the great pieces of mathematical insight pertaining to braids, the Nielsen-Thurston classification theorem~\cite{MR15791,MR956596,MR3053012,MR964685}.  With this tool we can precisely define what we mean by simple or complex braids, as well as any combination of these categories.  Complex braids, called pseudo-Anosov, have an important topological invariant, the braid dilation (or the log of this, the topological entropy~\cite{MR175106,MR274707,MR0255765}).  Roughly speaking, the topological entropy, which we introduce in section~\ref{Sec:Mixing}, captures the exponential rate of stretching of material curves due to the motion of the points forming the braid.  Since stretching and folding underlie all kinematic mixing mechanisms~\cite{aref_1984,ottino1989kinematics} in fluid dynamics, the topological entropy constitutes an important measure of mixing.

As the braid dilation is an algebraic number, many researchers have considered  the interesting question of what braids give minimal, though non-zero, topological entropy~\cite{MR2828128,MR2795241,10.2140/agt.2006.6.699}.  While the analogous maximization question is more important to applications, as it concerns maximizing mixing, it does require additional assumptions to make it well-posed.  We can increase the topological entropy of a braid without limit, simply by making the braid longer.  What is needed is a way to normalize the topological entropy to get a quantity that reflects the efficiency of topological entropy production, or the rate of mixing.  For braids that arise from dynamical systems, we can normalize by the time it takes for the trajectories to form the braid.  Alternatively, we could normalize by the path-length of trajectories, by the work done on the system, or to suit the optimization requirements of an engineering problem.  While physically meaningful, these normalizations are not intrinsic to the braid itself.  By rescaling time, e.g. by increasing particle speeds, we can change the time-normalized topological entropy without changing the underlying complexity of the braid.  We would like a normalization which is intrinsic to the braid (or at least the braid presentation).

A standard way to represent braids algebraically is to use Artin braid generators~\cite{MR19087}.  Given a set of points equally spaced along a line in the plane, each generator represents a clock-wise or counter clock-wise switch of two adjacent points.  A braid word is a sequence of these generators.  We can normalize the topological entropy of a braid word by the number of generators it contains, giving the topological entropy per generator (TEPG).  However, given that multiple generators can be executed simultaneously if no two of them are pair-wise adjacent (share a point), the TEPG doesn't really capture normalization by an intrinsic and discrete notion of time.  A better normalization starts with braid operations, which consist of sets of generators that can be executed simultaneously.  The topological entropy per operation (TEPO) follows from normalizing by the minimal number of operations needed to define a given braid word~\cite{MR2861264}.  We will try to find braids that maximize mixing efficiency as expressed by the TEPO in section~\ref{Sec:Search}.

Part of the appeal of the Artin presentation of braids lies in the relatively small number of braid relations needed to relate braids and braid words.  However, from the perspective of modeling point movement on the plane, there are some distinct disadvantages.  Chief among them is the fact that the canonical way of associating geometric movement in the plane with Artin generators requires that the points lie on a line.  The resulting movement does not fully take advantage of the two-dimensionality of the plane.  To remedy this, we consider a generalization of the Artin algebraic presentation using planar graphs embedded in orientable surfaces.  With this presentation of surface braids, introduced in section~\ref{Sec:Braids}, we can analogously define our measure of mixing efficiency, the TEPO.
 
Computing the topological entropy and TEPO for a given braid is a non-trivial problem, which has inspired an interesting history of algorithm development.  In section~\ref{Sec:MixEff} we layout the key ideas for our algorithm.  In short, we use a triangulation to create a coordinate system for closed curves, calculate the action of braid generators on these coordinates, and obtain the exponential increase in curve length which constitutes topological entropy.  This approach is very flexible, allowing for any combination of orientable surface and graph presentation.  It is also fast, enabling brute force searches for maximum TEPO braids.

As an example of how our algorithm works, we consider braids on planar lattice graphs in $\mathbb{R}^2$.  In particular, we look at minimal torus models of these lattice graphs.  These models, introduced in section~\ref{Sec:braidsasmodel}, constitute a more natural starting point for generating point motion throughout the plane as compared to Artin braids.  We walk though our algorithm for the case of the two point square lattice graph embedded in a torus in section~\ref{Sec:Example}.  

Furthermore, we show our search results for maximal TEPO braids on the six simplest torus models in section~\ref{Sec:Search}.  From this search we identify one simple braid which we conjecture to have the maximum TEPO for all braids on planar lattice graphs.  In section~\ref{Sec:MaxTEPOBraid}, we introduce various properties of this braid, including an analytical form for its TEPO, its measured invariant train tracks~\cite{MR1144770}, and its veering triangulation structure~\cite{Agol2011IdealTO}.

Finally, in section~\ref{Sec:Connections}, we come back to the connection with applications.  We see that the max TEPO braid forms the basis for an efficient stirring rod mixing mechanism,  and that it might be seen in the motion of topological defects in active nematics.  In appendix~\ref{Appendix:OtherLatticeGraphs}, we fill in some of the details of the other torus models included in the search, and share where to get the python implementation of the braiding algorithm for the six torus models.  In appendix~\ref{Appendix:gensq}, we introduce the algorithm for the general torus model of the square lattice graph.  We also highlight some of the difficulties, like finding perfect matchings on the graph, that arise when considering larger torus embeddings.

\section{Braids}
\label{Sec:Braids}

The classic braid group, $\mathbf{B}_n$, has many different realizations~\cite{MR0425944,kassel2008braid}: as the fundamental group of certain configuration spaces, as the mapping class group of the $n$-punctured disk, as homotopy classes of geometric braids, and as algebraic  braids with the Artin presentation.  Here we consider a generalization of geometric braids - surface braids\cite{Bellingeri2001OnPO}.  In particular, we are interested in  algebraic presentations of surface braids based on surface-embedded graphs and Artin-like generators.

Consider the orientable surface, $S$, and a set of $n$ distinguishable points $\mathcal{P} = \{ P_1,P_2, \cdots, P_n \}$ in $S$.  A geometric surface braid on $S$ is defined by a set of paths $\Gamma = (\gamma_1,\gamma_2, \cdots, \gamma_n)$, $\gamma_i\colon [0,1] \to S$, such that for $i,j \in \{1,2,\cdots,n\}$, $\gamma_i(0) = P_i$, $\gamma_i(1) \in \mathcal{P}$, and $\gamma_i(t) \neq \gamma_j(t)$ for all $i \neq j$ and all $t \in [0,1]$.  When $S = \mathbb{R}^2$, we recover the usual geometric braids.  Here we can treat $t$ as an orthogonal spatial dimension, where the paths $\gamma_i$ are geometric strands, just like strands of hair entwining to form what most non-mathematicians associate with the word ``braid".  Alternatively, $t$ can be a time parameter, where the paths now describe the movement of points on the given surface.  We will take this point of view, which directly connects braids and dynamics.

Since we are interested in the topological features of how points wind about each-other, and not the specific geometry of their paths, $\gamma_i$, we will focus on surface braid groups.  To move from geometric surface braids to surface braid groups, we consider homotopy classes of geometric surface braids, where any two geometric surface braids in a class can be connected by a one-parameter family of geometric surface braids.  In other words, a homotopy class represents all collections of paths, $\Gamma$, that can be continuously deformed into one another without path intersections. This, along with the natural notion of a product from path concatenation, defines a group structure.  The resultant surface braid groups, $\mathbf{B}(n,S)$ (e.g. $\mathbf{B}_n = \mathbf{B}(n,\mathbb{R}^2)$), do not depend on the specific set of points, $\mathcal{P}$.  We will simply refer to an element of a surface braid group, i.e. a homotopy class, as a braid.

Surface braid groups capture the topological aspects of point motion that we are interested in, but are not easy to work with computationally.  To help with this, we will develop an algebraic presentation of surface braid groups that is amenable to computational manipulations.  We start with the standard Artin algebraic presentation~\cite{MR19087}, which gives $\mathbf{B}_n$ in terms of $n-1$ braid generators $\sigma_1,\cdots,\sigma_{n-1}$ (as well as their inverses, $\sigma_i^{-1}$), and a small set of ``braid relations": $\sigma_i\sigma_i^{-1} = \mathbf{1}$, $\sigma_i\sigma_j = \sigma_j\sigma_i$ for $|i-j|\geq 2$, and $\sigma_i\sigma_{i+1}\sigma_i = \sigma_{i+1}\sigma_i\sigma_{i+1}$.  A braid word is an ordered sequence of braid generators, $\beta = \sigma_{i_1}^{\epsilon_1}\sigma_{i_2}^{\epsilon_2}\cdots\sigma_{i_k}^{\epsilon_k}$, where $i_j  \in \{1,2, \cdots, n-1\}$ and $\epsilon_j \in \{+1,-1\}$.  Two braid words represent the same braid if there exists a sequence of braid relations that connect the two.  This braid word problem has been solved in many different ways~\cite{MR2462117}, though its solution does not concern us here.

Before introducing the graph-based algebraic presentation for surface braids, we point out that the group isomorphism between geometric braid groups and Artin algebraic braid groups mirrors the ways in which braids have generally been applied to dynamical systems.  Associating an Artin algebraic braid word with a given geometric braid involves projecting the strand positions onto an axis and encoding changes in the strand ordering, noting how strands cross, with a sequence of Artin generators.  This is precisely the process used to go from trajectory data passively generated by a dynamical system to an Artin braid word as a discrete encoding of the flow.   In the opposite direction, there is a canonical way to assign a geometric braid to each Artin algebraic braid word.  Consider $n$ points equally spaced out on the x-axis in $\mathbb{R}^2$, starting on the left at the origin.  For each generator $\sigma_i$, we rigidly rotate the $i^{th}$ pair of points by half a turn about their geometric center in a counter clockwise (CCW) manner (or clockwise for $\sigma_i^{-1}$).  These half Dehn twists, ordered in time as they appear left to right in the braid word, generate a set of point trajectories.  In this way we can impose geometric motion, of e.g. mixing rods in a fluid, that is consistent with an algebraic braid word.

However, these two constructions pose some problems for using braids to model or generate motion in 2D.  First of all, the linear projection step results in Artin generators that do not model motion locally:  Two points that switch positions along the projection axis might actually be on opposite sides of the system.  Furthermore, the Artin representation of two points that switch positions locally depends on the relative location of other points in the system.  In the other direction, going from algebraic braid words to geometric braids, the restriction to colinear points precludes generating patterns of movement that are fully two dimensional.

\begin{figure}[htbp]
		\center
		\includegraphics[width = 1.0 \linewidth]{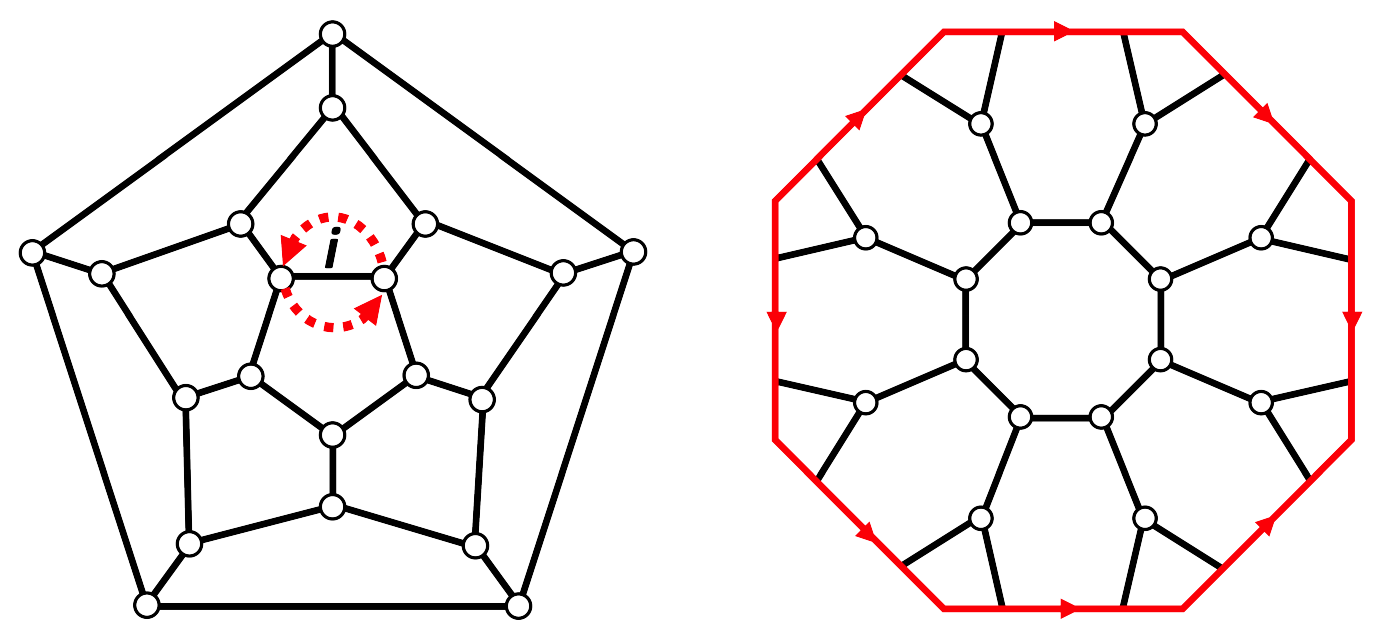}
		\caption{Two examples of graphs embedded in surfaces.  On the left is the skeleton graph of the dodecahedron. Though shown as a planar graph, this can also be embedded on the sphere (considering the rest of the plane as a pentagonal face through one point compactification).  On the right is a graph on the two-torus with octagonal faces (identify opposite sides of the red octagon to get the genus 2 surface).  Both graphs are regular maps.  On the left, edge $i$ is labeled, and the point motions associated with the generator $\sigma_i$ are indicated by the red arrows.}
		\label{Fig:GraphEmbed}
	\end{figure}

We would like a braid presentation whose  generators are local, and which could encode local point motions.  As a first step, consider the pair-wise rotations of points on a line associated with Artin generators.  We can conceive of these points on a line as a graph, where the points are vertices and the sections of the line between points are the edges.  Now each edge in the graph is associated with a generator (and its inverse); the graph encodes which pairs of points may swap places in one move.  As a natural extension of this idea, we consider analogously defined braids using more complex graphs\cite{Sergiescu1993}.  For example, an annular braid is defined on a graph similar to the one above, but with one extra edge connecting the two end points.  These graphs need not lie on the plane, but can be embedded in various surfaces\cite{mohar2001graphs}, see fig.~\ref{Fig:GraphEmbed} for examples.  

More formally, consider the graph $\mathcal{G} = (\mathcal{V},\mathcal{E})$, where $\mathcal{V}$ is the vertex set ($|\mathcal{V}| =  n$), and $\mathcal{E}$ is the edge set.  An embedding of $\mathcal{G}$ in an orientable surface $S$, associates each vertex, $v_i$, with a point $P_i \in S$ and each edge $e_j$ with a curve $\eta_j \colon [0,1] \to S$, such that for each edge $e_j$ with adjacent vertices, $v_a$ and $v_b$, we have $\{\eta_j(0),\eta_j(1)\} = \{P_a,P_b\}$.  We consider graphs that are planar relative to $S$, i.e. there exists an embedding such that $\eta_i(t) \neq \eta_j(s)$ for $i \neq j$ and all $0 < t < 1$, $0 < s < 1$.  Again, we are not interested in the geometric particulars of the embeddings, and so consider all homotopically equivalent embeddings to be the same.  For brevity, we will simply refer to a homotopy class of surface-embedded planar graphs as just graphs.  

In this context, a braid generator, $\sigma_i$ (or $\sigma^{-1}_i$), is defined to be a CCW (resp. CW) switch of the two points, $P_a$ and $P_b$, whose associated graph vertices, $v_a$ and $v_b$, are adjacent to the graph edge $e_i$; e.g. see the left side of fig.~\ref{Fig:GraphEmbed}.  Again, a braid word is an ordered sequence of these braid generators.  For this to constitute an algebraic presentation of the surface braid group $\mathbf{B}(n,S)$, we also require a set of braid relations.  There is such a set\cite{Sergiescu1993}, but for our purposes we need only the identity relation, $\sigma_i\sigma_i^{-1} = \mathbf{1}$, and the commutation relation,  $\sigma_i\sigma_j = \sigma_j\sigma_i$ if $\mathcal{V}(e_i) \cap \mathcal{V}(e_j) = \emptyset$, where ${\mathcal{V}(e_k)}$ is the set of vertices adjacent to edge $e_k$.

This presentation is much more flexible than Artin braids for modeling point motion on various surfaces.  However, there is a trade-off, as there no longer exists a general canonical procedure for passing from geometric braids to algebraic braid words, like there is with the projection step for Artin braids.  The opposite procedure is possible, and we can always go from a braid word to a movement of points that geometrically realizes this braid on the embedded graph.  In light of this, we mainly treat braids as a discrete model of the possibilities for point motion in the plane.

Finally, braid operations are defined as sets of generators that pairwise commute, and therefore can be executed simultaneously.  We consider maximal braid operations, where as many vertices as possible participate in swaps.  Enumerating the number of possible operations for a given graph is an interesting combinatorial problem (related to matchings in graph theory), which is explored further in appendix~\ref{Appendix:gensq}.  Each braid word can be written as a time-ordered sequence of braid operations.  The minimum number of operations needed to express a braid word constitutes a discrete notion of the time it takes to execute the braid word.  Indeed, this ``topological time" should be a braid invariant if, for a given surface braid, we minimize over all compatible braid words in every graph presentation with the same number of vertices.  We will be using the number of braid operations to normalize the topological entropy.

\section{Lattice Graphs and Torus Models}
\label{Sec:braidsasmodel}

In section~\ref{Sec:MixEff} we introduce an algorithm to compute the topological entropy of surface braids expressed using graph generators.  While this method applies generally, we would like to concretely illustrate how it works using some specific examples that are both intrinsically interesting and important for applications.  In this section we introduce and motivate these examples.

We will focus on $\mathbb{R}^2$ as an example surface, since the majority of braid applications mentioned in section~\ref{Sec:Intro} involve the movement of points on the plane.  Though later we will use the torus to model spatially periodic motion on the plane.  For the sake of simplicity, and to allow for motion that is as isotropic and homogeneous as possible, we require the graphs to be maximally symmetric.  More specifically, we choose graphs that are arc-transitive, also called symmetric graphs (or regular maps if referring to the embedding).  A graph is symmetric, if for every two pairs of adjacent vertices, say $v_a-v_b$ and $v_c-v_d$, there is a graph automorphism, $\Psi$, such that $\Psi(v_a) = v_c$ and $\Psi(v_b) = v_d$.  Thus, we can map every edge to every other edge two ways, while mapping the graph back onto itself.  Fig.~\ref{Fig:GraphEmbed} shows two example symmetric graphs, though they are embedded on the sphere and two torus respectively.   For the moment, we consider infinite graphs, acknowledging that these are computationally infeasible and that we will later need to introduce finite models for them. Finally, we require that the geometric embedding itself is just as symmetric as the graph.  That is to say, each graph automorphism can be realized on the embedded graph through rigid rotations and translations of the plane.  As a consequence, every edge has the same length, and the generators will model points moving with constant speed. 

\begin{figure}[htbp]
		\center
		\includegraphics[width = 1.0 \linewidth]{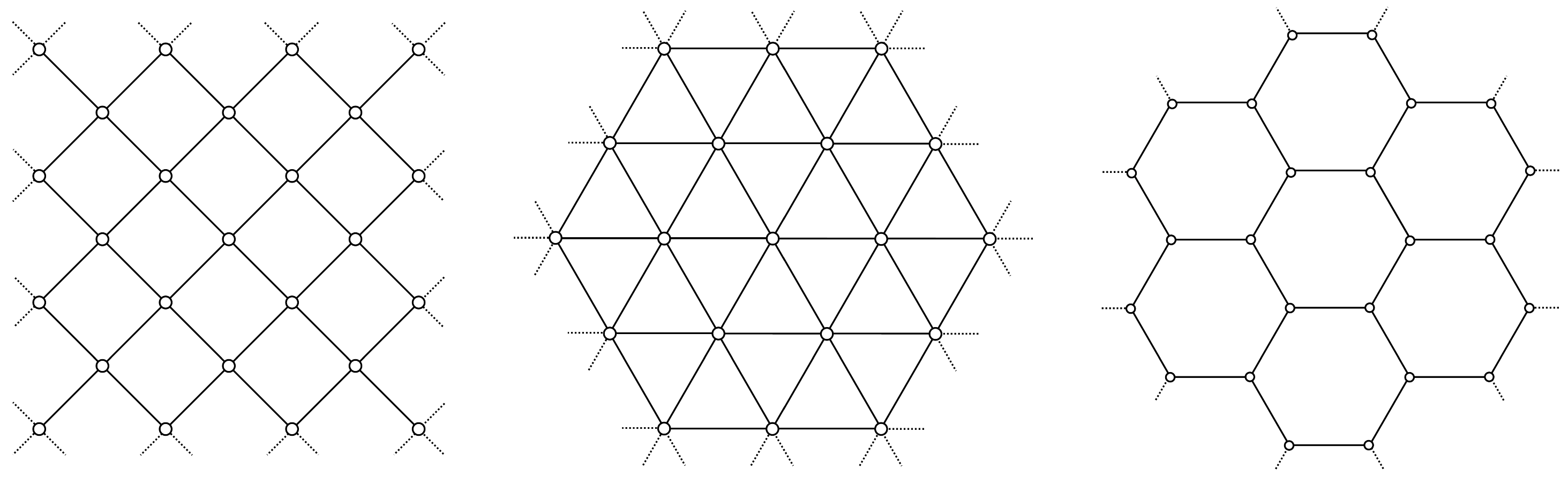}
		\caption{The three lattices (square, triangular, and hexagonal), which we will consider.  These constitute the most symmetric graphs on the plane.}
		\label{Fig:LatticeTypes}
	\end{figure}

In short, we consider the simplest, most symmetric, infinite graphs embedded in the plane.  We will refer to these as lattice graphs, which come in three types: square, triangular, and hexagonal, as seen in fig.~\ref{Fig:LatticeTypes}.

	\begin{figure*}[htbp]
		\center
		\includegraphics[width = 1.0 \linewidth]{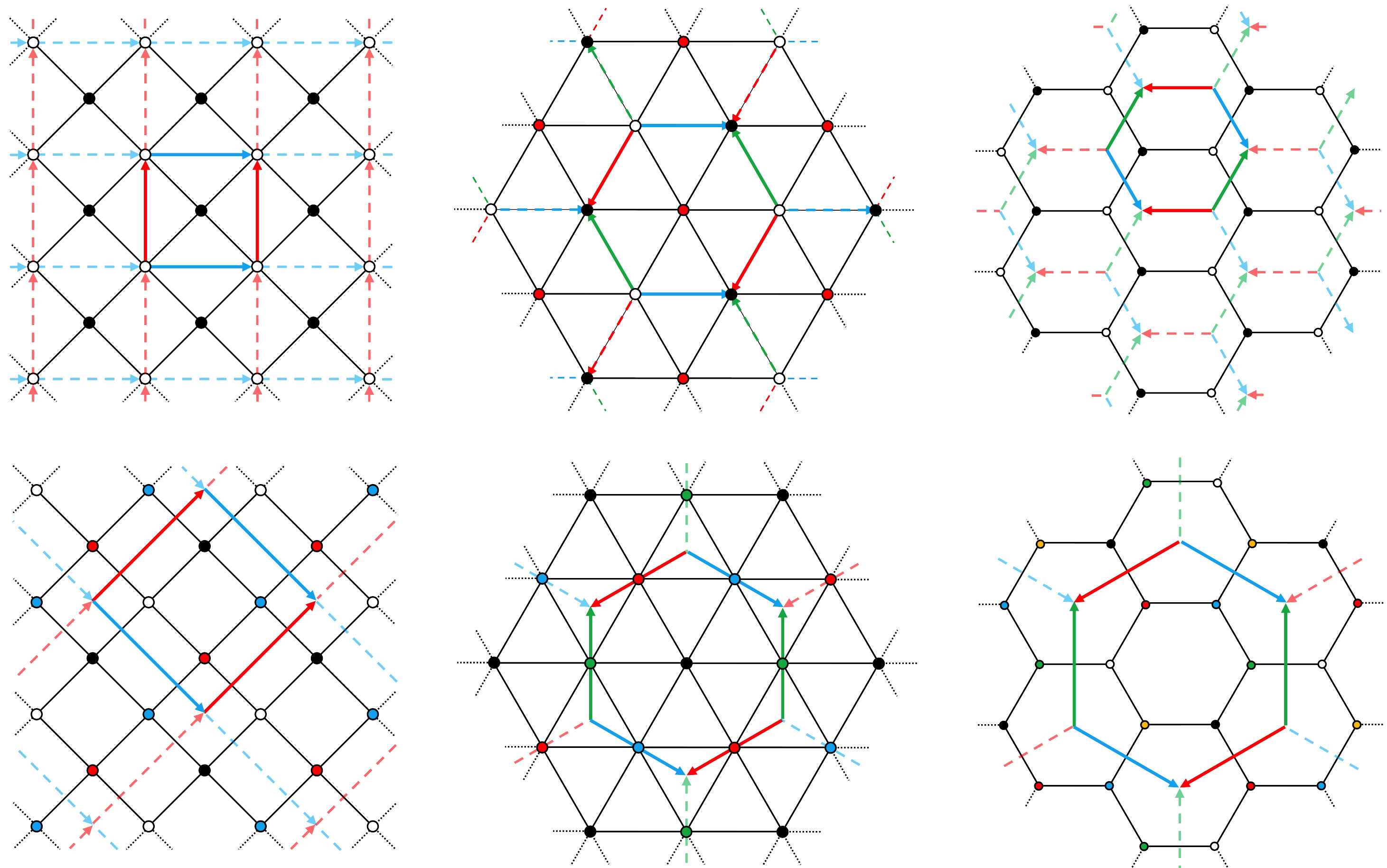}
		\caption{Torus embedding models for lattice graphs.  The smallest two dynamically non-trivial ways to embed square, triangular, and hexagonal lattice graphs on the torus are shown.  Square fundamental domains are shown with blue and red sides, while hexagonal fundamental domains have an additional green side (for each, identify like sides to get the torus).}
		\label{Fig:FundamentalDomains}
	\end{figure*}

Braids on lattice graphs have a countably infinite number of braid generators, and so any practical attempt to write down a braid must necessarily consider spatially periodic braids.  For this, we consider a fundamental domain, on which the finite number of braid generators are defined, and each translated copy of this domain used to tile over the entire embedded lattice graph replicates these generators.  Since we are tiling the plane with translations only, we can only use square and hexagonal fundamental domains.  Furthermore, there are a discrete number of ways to pick a fundamental domain such that the lattice graph structure within the fundamental domain is itself preserved upon translation to every tile in the covering.  These are distinguished by the number of vertices contained in the fundamental domain.  For instance, for the square lattice graph, the two smallest (dynamically interesting) fundamental domains are themselves squares, and enclose 2 and 4 points, see fig.~\ref{Fig:FundamentalDomains}.  Both the triangular and hexagonal lattice graphs require a hexagonal fundamental domain, and the two smallest fundamental domains for each contain 3 and 4 points, and 2 and 6 points respectively, see fig.~\ref{Fig:FundamentalDomains}.

Since the periodicity of the tilings allow us to identify the opposite sides of the square fundamental domain and separately the hexagon fundamental domain, the natural space that these periodically restricted lattice graph braids live on is the torus.  Indeed, the six examples in fig.~\ref{Fig:FundamentalDomains} are well known from the study of regular maps on the torus~\cite{RegularMapsSurvey}.  From this viewpoint, we are simply considering braids defined on the two smallest embeddings of the three different lattice graphs in the torus.  We will refer to these as braids on torus lattice graphs, or more succinctly lattice braids.  Others  have previously investigated the topological entropy of torus braids~\cite{MR2299637}, though they used a different braid presentation and different method for finding topological entropy.

Now, for each choice of lattice graph type and embedding type, we have a small set of edges in the fundamental domain, and therefore a small set of braid generators.  Using the graphs structure to determine generator adjacency, we can enumerate the braid operations.  As a reminder, braid operations are matchings of the graph (with the additional information about CW/CCW swaps), collections of edges that do not share any vertices.  While any such collection of generators constitutes an operation, we will only consider operations that correspond to perfect matchings, where each point is matched.  The number of possible braid operations draws on the combinatorics of perfect matchings.  As an example, consider the square lattice graph.  For the torus embedding with two points, there are four edges, and therefore four generators (as well as the four inverses).  Since a single generator involves both of the points, it is also an operation, and there are 8 possible braid operations.  For the 4 point torus embedding, there are 8 edges/generators.  There are 8 distinct matchings, which, with 4 possible CCW/CW combinations for the two constituent generators, gives a total of 32 possible braid operations.  These operations form the alphabet which we use to build braid words.

\section{Mixing}
\label{Sec:Mixing}

As we have seen, the topologically relevant features of the motion of points in the plane are represented algebraically by braids.  This motion stretches and folds the the surrounding medium in which the points are embedded, and the amount of mixing can be gleaned from attributes of the braid.  In particular, we will be able to associate with each braid a topological invariant, the topological entropy~\cite{MR175106,MR274707,MR0255765}, that serves as a measure of mixing.

As a starting point, we go back to Artin braids and their realization as the mapping class group~\cite{MR2850125} of the punctured disk.  The mapping class group of a surface consists of elements, each corresponding to a class of homeomorphisms of the surface to itself (potentially permuting any punctures), that are equivalent under homotopy.  The main tool for making sense of mapping class groups is the Nielsen-Thurston classification theorem~\cite{MR15791,MR956596,MR3053012,MR964685}.  Roughly, this says that for each element of the mapping class group, we can choose a representative map (homeomorphism) with well defined properties.  These representative maps are either finite order, pseudo-Anosov (pA), or reducible.  We will only be concerned with pA maps, though briefly, finite-order maps are those for which some power of the representative map is the identity, and reducible maps can be cut up into combinations of finite order and pA maps.

Pseudo-Anosov maps have particularly rich behavior, and the braids associated with them (for the punctured disk) have good mixing properties.  For the representative pA map, there exists two measured foliations~\cite{MR3053012}, transverse to one-another, and a positive real number - $\lambda$ - the dilation.  Under the action of the map, the unstable foliation's transverse measure is multiplied by $\lambda$, and that of the stable foliation is multiplied by $1/\lambda$.

\begin{figure}[htbp]
		\center
		\includegraphics[width = 1.0 \linewidth]{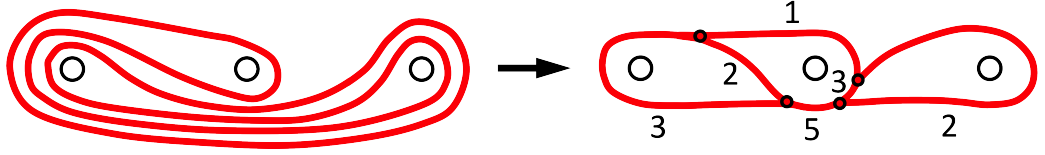}
		\caption{A closed curve stretched about three points can be represented by a weighted train track.  Here the measure on each train track branch encodes the number of parallel sections of the curve that have been bunched together.}
		\label{Fig:TrainTrackEx}
	\end{figure}

We will recast this property of pA maps in terms of measured train-tracks (instead of foliations), which will underpin our algorithm for calculating the mixing efficiency for braids.  For our purposes, a train track~\cite{MR1144770}, $\tau$, consists of a geometrically realized graph, (typically vertices are called switches, $v_i  \in \mathcal{V}(\tau)$, and edges are called branches, $b_i \in \mathcal{E}(\tau)$), with the condition that three branches meet at each switch and are tangent (two branches on one side of the switch and one branch on the other), see fig.~\ref{Fig:TrainTrackEx}.  Measured train tracks assign to each branch, $b_i$, a real number, its weight - $\mu(b_i)$, such that at each switch the weights are conserved (``switch  conditions" - weights of the two branches on one side combine to give the weight of the branch on the other side, or $\mu(b_1)+\mu(b_2) = \mu(b_3)$).  Train tracks were created to represent non-overlapping sets of closed curves on surfaces, where parallel sections of curves are bunched together into a train track branch with the edge weight representing the number of curves that were bunched together, fig.~\ref{Fig:TrainTrackEx}.  Measures with non-negative real valued weights can also be used, provided they satisfy the switch conditions.

\begin{figure}[htbp]
		\center
		\includegraphics[width = 1.0 \linewidth]{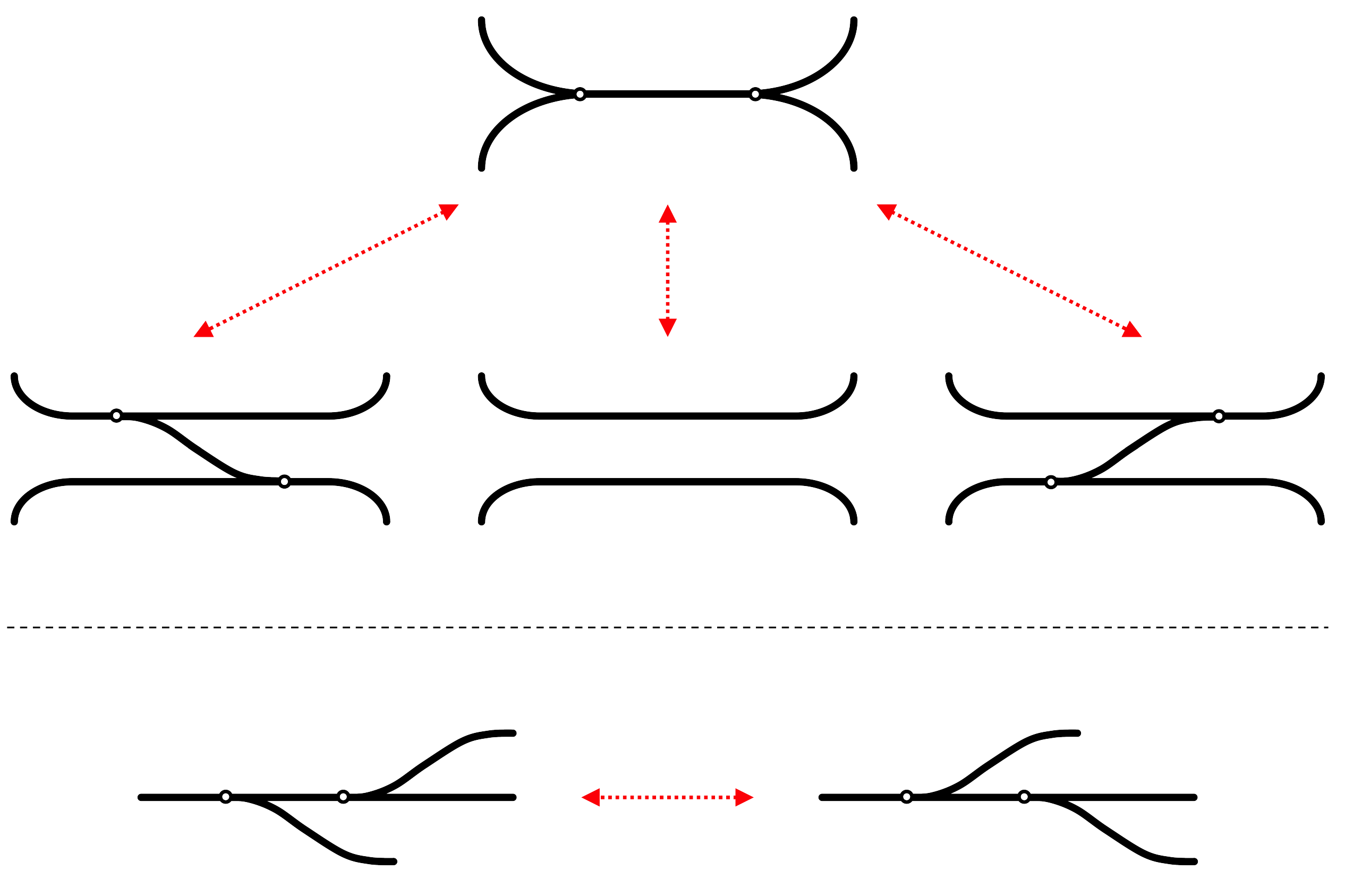}
		\caption{Train track operations.  Top: Splitting a branch results in one of three topologically distinct configurations.  The choice of configuration depends on the measure.  In contrast, the train track changes due to folding branches (the inverse of splitting) are independent of the measure.  Bottom: A slide move and its inverse (also a slide move).}
		\label{Fig:TrainTrackMoves}
	\end{figure}

The main operations on train tracks are splits, folds, and slides, as can be seen in fig.~\ref{Fig:TrainTrackMoves}.  Folds and splits are the inverse of one another, while a slide can be undone by another slide.  Slides and folds do not require any knowledge of the branch weights, while choosing between the three topologically distinct splitting options requires the measure.

For a representative pA map, there exists a measured train track such that a series of folding, splitting, and slide operations, as well as movements, map the train track graph back onto itself, and the weight of each brach has increased by a multiplicative factor of $\lambda$, the braid dilation.  This is a combinatorially efficient way of representing the unstable measured foliation; there is likewise a train track representing the stable measured foliation.

This property of the pA representative map is important, as it forces at least as much stretching as the dilation in every single map in this class.  This means that no matter the material or underlying dynamics of the system with moving points, as long as the trajectories of the points form a pA braid, there must be some material curves that stretch out by a factor at least as large as the braid dilation.

Our measure of mixing, the topological entropy, is very closely related to the dilation.  For a homeomorphism, the topological entropy is, loosely, the exponential rate of increase in distinguishable orbits with mapping iteration~\cite{MR0255765}.  For a braid, the topological entropy is simply $h = \log(\lambda)$.  More intuitively, it is the asymptotic exponential stretching rate of material lines that non-trivially enclose the points which correspond to the braid strands.  This suggests a practical method for finding the topological entropy of a braid $\beta$, assuming that we know how to calculate the action of this braid on a train track, as well as reliably encode the train track information.  We start with a measured train track, $\tau_0$, which represents some simple essential closed curve.  The action of the braid results in a new measured train track, $\tau_1 = \beta(\tau_0)$, and by repeated application of the braid we get a sequence of train tracks ($\tau_0,\tau_1, \cdots$).  Say $W_k = \Sigma_{b_i \in \tau_k} \mu(b_i)$ is the sum of the weights over every branch of $\tau_k$.  Then 
\begin{equation}
 h = \lim_{k \rightarrow \infty} \log(\frac{W_k}{W_{k-1}})
\end{equation}
Practically speaking, this will converge within some tolerance after a relative small number of steps for pA braids.  In the next section we will see how to represent the train track and how to calculate the action of a braid.

Now, the topological entropy is extensive in that $h(\beta^n) = n \times h(\beta)$, and therefore is not a good measure of mixing efficiency.  Executing a braid twice will double the topological entropy and will double the number of moves needed to execute this braid.  A measure of mixing efficiency which is intensive, i.e. $\overline{h}(\beta^n) = \overline{h}(\beta)$, is given by normalizing the topological entropy by the minimum number of braid operations, $|| \beta ||$, needed to write the braid word.  As a reminder, braid operations consist of sets of pairwise commuting braid generators.  The topological entropy per operation (TEPO = $\overline{h}$) is our measure of mixing efficiency.

For the Artin braid group $\mathbf{B}_3$ (the smallest braid group with pA braids), the braid with maximal TEPO~\cite{MR2861264} is $\beta_{\phi}  = \sigma_1\sigma_2^{-1}$.  It has $\overline{h}(\beta_{\phi}) = log(\phi)$, where $\phi = \frac{1+\sqrt{5}}{2}$ is the golden ratio.  Similarly, for braids on an annulus with an even number of points, $2n$, the max TEPO braid~\cite{MR2861264} is $\beta_{\delta_s} = \sigma_{odd}\sigma_{even}^{-1}$, where $ \sigma_{odd} =  \sigma_{1}\sigma_{3}\cdots\sigma_{2n-1}$ and $ \sigma_{even} =  \sigma_{2}\sigma_{4}\cdots\sigma_{2n}$ are braid operations.   Here the value of the topological entropy per operation is $\overline{h}(\beta_{\delta_s}) = log(\delta_s)$, where $\delta_s = 1+\sqrt{2}$ is the silver ratio.  We will find that our candidate maximal TEPO lattice graph braid also has a very pleasing analytical form for its mixing efficiency.

\section{Calculating Mixing Efficiency}
\label{Sec:MixEff}

In order to calculate the topological entropy per operation for an arbitrary braid, we must establish a way to represent measured train tracks and the rules for updating train tracks due to the action of braid generators.  In this section we give an overview of a general algorithm to achieve this for surface braids with braid generators defined by embedded graphs.  In the subsequent section we will work through this method using a concrete example.

The classical computational approach to this problem involves the Bestvina-Handel algorithm~\cite{MR1308491}, which is itself an algorithmic proof of the Nielsen-Thurston classification theorem~\cite{MR15791,MR956596,MR3053012,MR964685} of mapping class groups.  This algorithm gives maximal information (the invariant train track, the dilation, as well as  reduction curves for the reducible case), but is very slow, making the analysis of braids with even moderate numbers of strands impractical.  A much faster approach uses a special coordinate system, Dynnikov coordinates~\cite{MR1918864,MR2512607,MR2381961}, to encode loops in the $n$ punctured disk.  Braid generators then have fairly straight-forward update rules for their action on the Dynnikov coordinates.   While fast, this method~\cite{MR3456018} applies to the Artin presentation of braids on the disk, and does not directly apply to the surface braids that we are interested in.  Another approach, the E-tec algorithm~\cite{E-tec}, creates a dynamic triangulation of ensembles of moving points, and represents loops by edge weights on this triangulation.  E-tec is more suited to the analysis of geometric trajectory data, and not movement defined by discrete generators.  Our approach here takes ideas from Dynnikov coordinates and E-tec, which results in a fast, elegant algorithm that can be applied to surface braids with graph generators.

The starting point for our algorithm is a coordinate system for encoding an arbitrary measured train track, $\tau$, which in turn represents a set of non-overlapping closed curves.   This is achieved by promoting the embedded graph, $\mathcal{G}$, to a triangulation of the surface, $\mathcal{T}$, by adding extra edges.  In the case of the triangular lattice graphs, the graph is already a triangulation, but for the square and hexagonal lattice graphs, extra edges are needed.  The exact way in which this is achieved is not very important, though triangulations that retain as many of the symmetries of the lattice graph itself are easiest to work with.  

Once we have a triangulation, we can associate with each edge in the triangulation, $e_i \in \mathcal{E}(\mathcal{T})$, an intersection coordinate, $E_i$.  This coordinate is given as $E_i = \min_{\tau' \sim \tau} \Sigma_{b_k \in I(e_i,\tau')} \mu(b_k)  $, where $ I(e_i,\tau')$ is the set of branches in $\tau'$ that transversely intersect the edge $e_i$.   Here, the minimum is taken over all train tracks, $\tau'$, which are equivalent to $\tau$ after a finite number of  folding, splitting, and slide operations as well as homotopies where the triangulation vertices act as obstructions.   In terms of closed curves, this simply states that if we pull the curves taught, so that they have minimum intersections with each triangulation edge, then the intersection coordinate for a triangulation edge counts the number of curves that intersect this edge.  

Now we can represent the measured train track, $\tau$, with the coordinate system $\vec{E} = (E_1,E_2, \cdots, E_{|\mathcal{E}(\mathcal{T})|})$.  Note that the coordinates can take values from the positive integers (to represent sets of closed curves), as well as from the non-negative real numbers (to represent asymptotic projective measured train tracks).  Also, note that, just as possible measures on the train track are restricted by the switch conditions, so too are the intersection coordinates restricted to those that satisfy triangle inequalities, $E_a \leq E_b +  E_c$, for the coordinates of each edge in a triangle (and all three inequalities for each triangle).  This means that the cardinality of these intersection coordinates is not minimal, and that the train track could be represented by a smaller set of coordinates.  Indeed, for a particular triangulation of the punctured disk, just such a reduction leads to Dynnikov coordinates.  We are trading a minimal representation for a simpler representation, which is easier to apply to general surfaces.  For a closed surface of genus $g$, with $n$ points, there will be $|\mathcal{E}(\mathcal{T})|  =  3(n+2g - 2)$ intersection coordinates.

	\begin{figure}[htbp]
		\centering
		\includegraphics[width = 0.6 \linewidth]{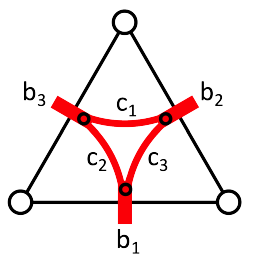}
		\caption{Constructing the canonical train track from intersection coordinates.  If the intersection coordinates obey the triangle inequalities, then we can calculate how each incoming train track branch splits and specify the measure for each connecting branch, eq.~\ref{Eq:TTtriangle}.}
		\label{Fig:CanonicalTT}
	\end{figure}

While many possible measured train tracks correspond to the same set of intersection coordinates, we can construct one canonical representative measured train track directly from the intersection coordinates.  Fig.~\ref{Fig:CanonicalTT} shows a generic triangle in the triangulation, with intersection coordinates $E_1,E_2,E_3$.  We draw one train track branch crossing each edge (labeled $b_1,b_2,b_3$, with measure $\mu(b_1) = E_1$, $\mu(b_2) = E_2$, $\mu(b_3) = E_3$), a switch on each branch interior to the triangle, and three new branches connecting the three switches (labeled $c_1,c_2,c_3$).  If these are valid coordinates, i.e. they obey the triangle inequalities, then we have
\begin{equation}
\begin{split}
\mu(c_1) &= (E_2+E_3-E_1)/2 \\
\mu(c_2) &= (E_1+E_3-E_2)/2 \\
\mu(c_3) &= (E_2+E_1-E_3)/2.
\end{split}
\label{Eq:TTtriangle}
\end{equation}
This canonical train track is unique for a given triangulation, and can be used as a representative of the class of measured train tracks related by train track operations and homotopies (relative to the points).  Thus, in order to describe the action of a braid generator on a measured train track class, it suffices to find this action on a canonical train track.  Furthermore, as a convention, we choose a reference triangulation so that we can compare canonical train tracks before and after the point movement associated with a braid generator.  Now, the generator's action can be captured by simple update rules for the reference triangulation edge weights.

	\begin{figure}[htbp]
		\centering
		\includegraphics[width = 1.0 \linewidth]{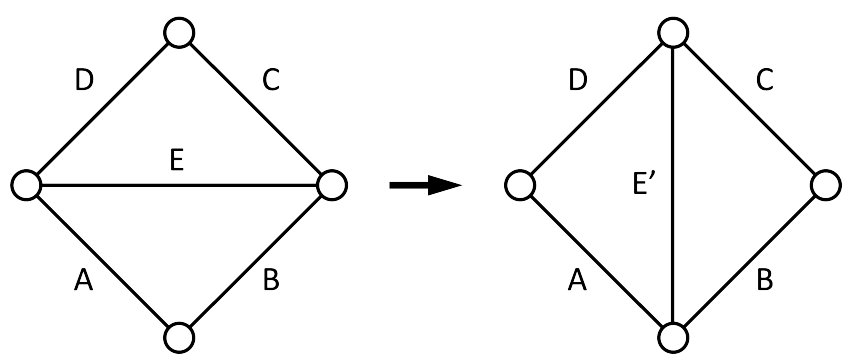}
		\caption{A triangulation flip.  The intersection coordinates A,B,C,D remain unchanged, while $E \rightarrow E^{'}$ according to eq.~\ref{Eq:Flip}.}
		\label{Fig:TriangleFlip}
	\end{figure}
	
The central building block used to construct these updates rules is a ``Flip" operation on the triangulation (sometimes called a Whitehead move).  In a flip we locally modify the triangulation by removing an edge shared between two triangles and replace it with an edge along the opposite diagonal of the resulting quadrilateral, see fig.~\ref{Fig:TriangleFlip}.  The flip induces a change in the intersection coordinate associated with the modified triangulation edge.  If the intersection coordinates of the edges surrounding the quadrilateral, in cyclic order, are $A, B,  C,  D$, and the initial diagonal edge has a weight of $E$ (as in fig.~\ref{Fig:TriangleFlip}), then the new edge's weight, $E'$ is
\begin{equation}
E' = \max(A+C,B+D) - E.
\label{Eq:Flip}
\end{equation}

	\begin{figure*}[htbp]
		\center
		\includegraphics[width = 1.0\linewidth]{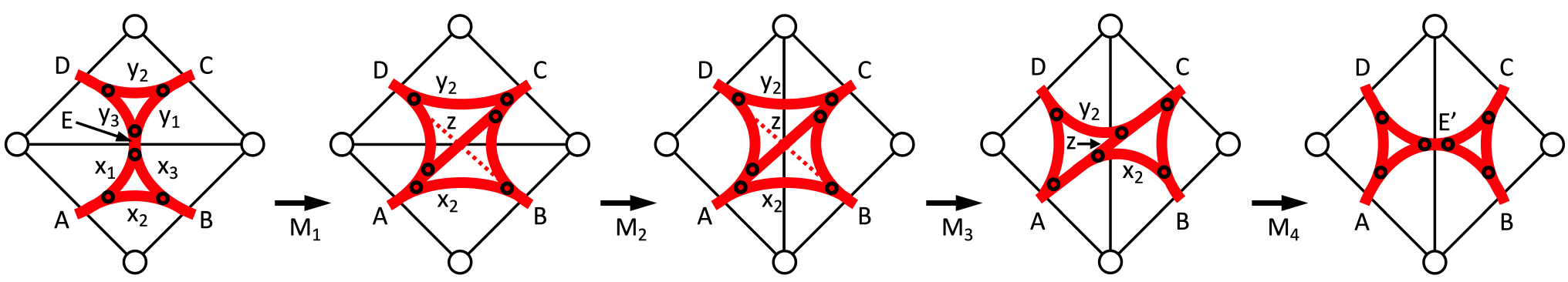}
		\caption{Derivation of the edge weight update rule, eq.~\ref{Eq:Flip}, for a flip operation, using a sequence of train track operations (splits, slides, and folds) on the canonical train track.  See eq.~\ref{Eq:FlipProof} and the accompanying discussion for details.}
		\label{Fig:FlipProof}
	\end{figure*}
	
\begin{figure}[htbp]
		\centering
		\includegraphics[width = 1.0 \linewidth]{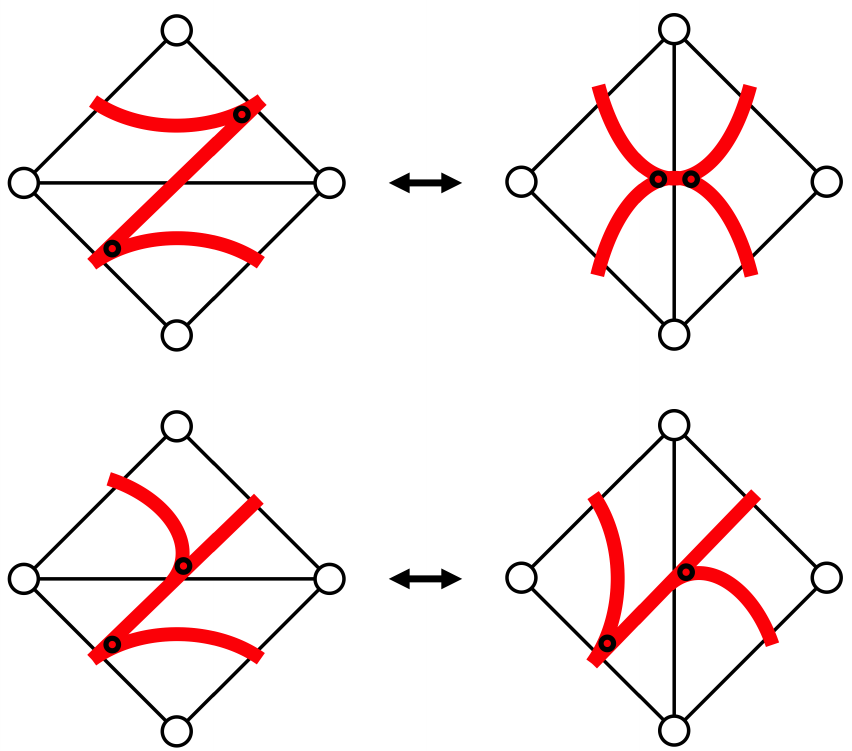}
		\caption{All train track operations are triangulation flips.  Top: train track folding operation (right arrow) or splitting operation (left arrow).  Bottom: a train track sliding operation. }
		\label{Fig:TTflip}
	\end{figure}

We derive this simple formula by relating the canonical train tracks before and after the flip using the train track moves shown in fig.~\ref{Fig:FlipProof}.  On the far left of this figure are two adjacent triangles, shown before the flip operation, and the local canonical measured train track corresponding to the edge weights $(A, B, C, D, E)$.  The measure on the train track branches interior to each triangle is given by the weights $(x_1,x_2,x_3)$ and $(y_1,y_2,y_3)$, which are themselves completely determined by the edge weights using eq.~\ref{Eq:TTtriangle}.  For move 1, $M_1$, we execute a split.  The new train track will have a branch, labeled $z$, which is oriented diagonally like ``/" (solid red line), oriented diagonally like ``\textbackslash" (dashed red line), or has a weight of zero (no line), depending on the weights.  A quick inspection tells us that $z = \max(y_3-x_1, y_1-x_3)$ encapsulates these three possibilities.  Next, $M_2$ simply flips the triangulation edge to its new, vertical, orientation.  In $M_3$, we execute two slide moves.  For clarity, we have dropped the dashed line and concentrated on only one of the split configurations from $M_1$ (the other two configurations are easy to draw).  Finally, $M_4$ executes a fold and gives us a measured train track that is in the canonical configuration for the new triangulation.  From inspection, $E' = z+y_2+x_2$.  To obtain eq.~\ref{Eq:Flip}, we execute a few simple algebraic manipulations using max-plus algebra~\cite{MR2681232,MR2188299} notation, where $\llbracket x+y \rrbracket \equiv \max(x,y)$ and $\llbracket xy \rrbracket \equiv x+y$ (with this notation we can use the normal distributive rules of algebra).	
\begin{equation}
\begin{split}
E' =&  z + x_2 + y_2  \\
 =& \max(y_3-x_1, y_1-x_3) + x_2 + y_2 \\
 =& \left \llbracket (\frac{y_3}{x_1}+\frac{y_1}{x_3})x_2 y_2 \right \rrbracket  \\ 
 =& \left \llbracket \frac{x_2 x_3 y_2 y_3 + x_1 x_2 y_1 y_2}{x_1 x_3} \right \rrbracket  \\ 
 =& \max(x_2 + x_3 + y_2 + y_3, x_1 + x_2 + y_1 + y_2) \\
 & - (x_1 + x_3) \\
 =& \max(B + D, A + C) - E.	
\end{split}
\label{Eq:FlipProof}
\end{equation}	
Here, the last step has used the train track switch conditions ($x_1+x_2 = A$, etc.).  Depending on the edge weights, this flip move can describe the standard folding, splitting, and sliding train track operations on the canonical train track (see fig.~\ref{Fig:TTflip}).
	
The flip operator has some interesting algebraic properties that warrant a few remarks.  First, the max operation keeps eq.~\ref{Eq:Flip} from being a linear equation, and therefore prevents the braid generator updating rules from being simple linear algebra.  However, for invariant train tracks, we can use the appropriate linear equation for each flip, and get an overall transition matrix for the action of the braid on the train track.  We will have more to say about this when we evaluate a particular braid in section~\ref{Sec:MaxTEPOBraid}.  Next, using the max-plus notation, the flip update rule can be written:
\begin{equation}
	\left \llbracket E^{'}E =  AC+BD \right \rrbracket
	\label{Eq:Flip2}
\end{equation}
Curiously, this equation has the same form as the Ptolomy relation for a quadrilateral.  Given any four points that lie on a circle forming a quadrilateral, the lengths of the two diagonals are related to the lengths of the four sides just as in eq.~\ref{Eq:Flip2}.  This, of course, is more than coincidence, and there are deep connections here that lead to tropical algebra, cluster algebras, and Teichmuller spaces~\cite{ClusterAlgebraDThurstonI,fomin2018cluster}.

\begin{figure}[htbp]
	\centering
	\includegraphics[width = 1.0 \linewidth]{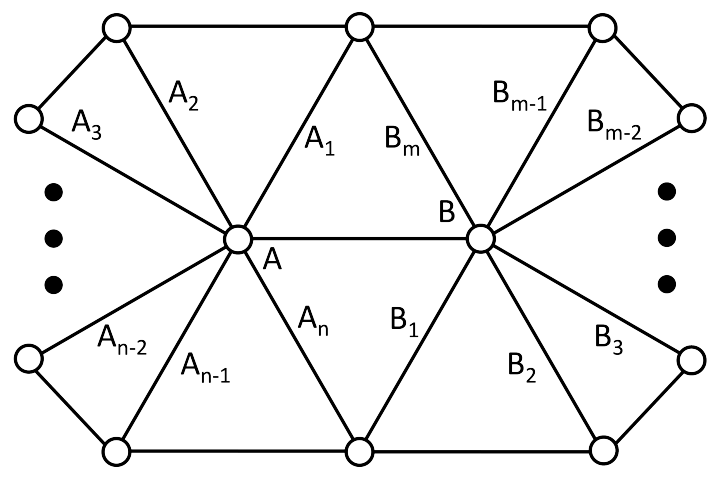}
	\caption{The triangles and edges adjacent to two points that participate in a swap.}
	\label{Fig:FlipPlanning}
\end{figure}

Now that we can represent a measured train track in canonical form (relative to a triangulation) using intersection coordinates, and can update these coordinates due to local retriangulation flips, we outline how to specify the action of a braid generator on triangulation coordinates.  In short we must produce a time ordered set of flip operations and point movements that map the original triangulation to itself and realize the point switching movement of the braid generator.  Any such set will work, but there is often a general procedure for creating this set.  Consider the set of triangles that are adjacent to either of the two points involved in the switch, as can be seen in fig.~\ref{Fig:FlipPlanning}.  If this set contains no duplicates (e.g., duplicates can happen for small triangulations on the torus due to the wrap-around boundaries), then we can implement the following flip sequence (for a CCW switch $\sigma$): $(A_1,A_2,\cdots,A_{n-1}, B_1, B_2,\cdots,B_{m-1})$.  Here, the label associated with an edge is used to represent a flip with this edge as the initial diagonal.  When this flip sequence is followed by a CCW swap of the points A and B, then the resultant triangulation will coincide with the original triangulation, and the associated intersection coordinates encode the action of this braid generator.  We will see an example in the next section where this procedure does not apply, yet we will still be able to come up with an appropriate flip sequence.

\begin{figure}[htbp]
	\center
	\includegraphics[width = 1.0\linewidth]{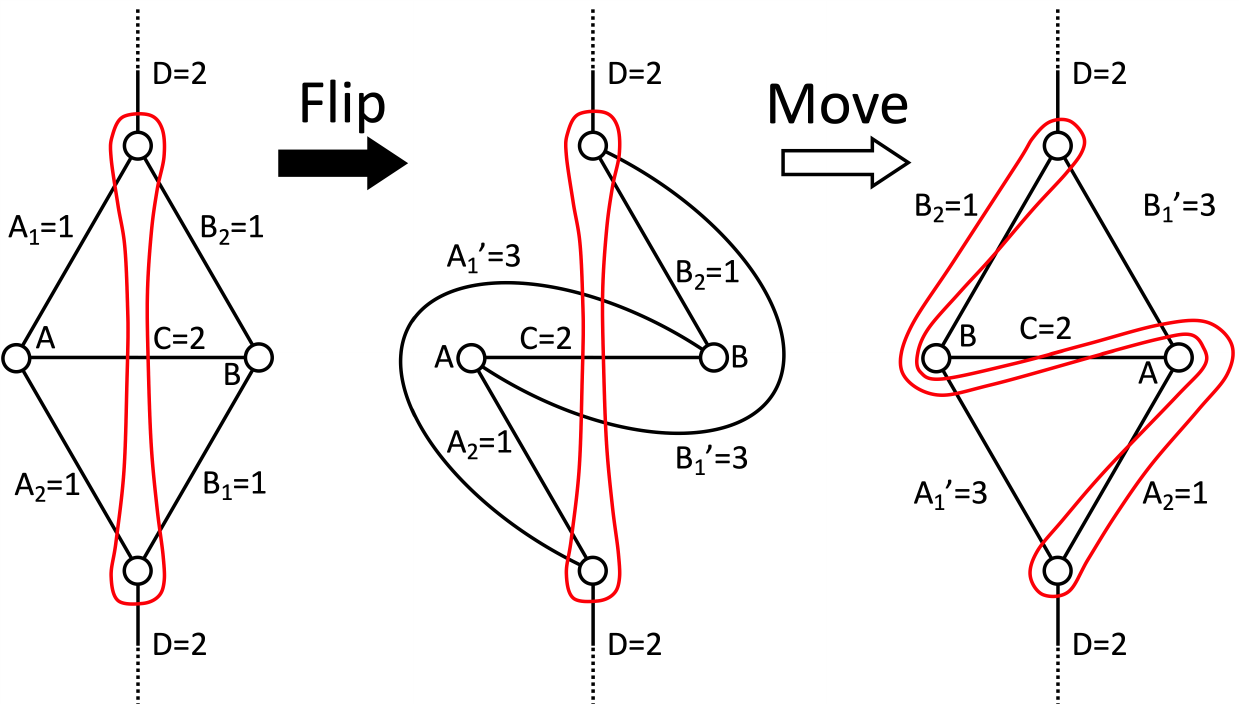}
	\caption{Simplest possible example of a flip sequence. Left: The black lines represent a tetrahedral triangulation of four points on the sphere (edge $D$ wraps around the other side), while the red loop represents an example closed curve (with indicated intersection coordinates).  Middle: Two flip operations have been performed, with the intersection coordinates of the new edges calculated using eq.~\ref{Eq:Flip}.  Right: The point movements corresponding to a braid generator (ccw swap of points A and B) are executed, and the final configuration of the example loop is shown.}
	\label{Fig:FlipSeqExSimple}
\end{figure}

For a concrete example of this procedure, consider fig.~\ref{Fig:FlipSeqExSimple}.  Here we have a tetrahedral triangulation of the sphere (edge $D$ wraps around the other side; $B_2-B_1-D$ and $A_2-A_1-D$ are both triangles).   Any measured train track can be represented in canonical form by the intersection coordinates $(A_1,A_2,B_1,B_2,C,D)$.  For the example loop shown in red, the coordinates are $(1,1,1,1,2,2)$.  There are two flips, $(A_1,B_1)$, in the flip sequence: flipping edge $A_1$ with quadrilateral edges $C,B_2,D,A_2$ and flipping edge $A_2$ with quadrilateral edges $C,A_2,D,B_2$.  These flips can be executed simultaneously and result in the middle triangulation.  The new edges have intersection coordinates which can be found using eq.~\ref{Eq:Flip}:  $A_1^{'} = \max(C+D,B_2+A_2) - A_1$ and $B_1^{'} = \max(C+D,B_2+A_2) - B_1$.  For the example loop we have: $A_1^{'} = \max(2+2,1+1) - 1 = 3$ and similarly $B_1^{'} = 3$.  In the final step, we move the points A and B counter clockwise about each other to realize the desired braid generator.  Comparing the intersection coordinates before and after these actions gives the update rules: $(A_1,A_2,B_1,B_2,C,D) \rightarrow (B_2,A_1^{'},A_2,B_1^{'},C,D)$, with $A_1^{'}$ and $B_1^{'}$ defined above.  Note that the final configuration of the example loop is consistent with the calculated new intersection coordinates.

For symmetric graphs (i.e. regular maps), it is sufficient to find the action of only one braid generator on the intersection coordinates.   The action of all other braid generators will be given by conjugating the action of this example generator with the action of various symmetries, such as translations, rotations, and mirror inversions.  For less symmetric graphs, the action of more braid generators will have to be directly calculated.  Once we have the action of each braid generator, we can build up the action of each braid operation, and from there, the action of an entire braid word.

\section{Example: 2 Point Square Lattice Graph on a Torus}
\label{Sec:Example}

%	\begin{figure}[htbp]
%		\centering
%		\includegraphics[width = 1.0 \linewidth]{Figures/SqLattice2ptTorus.pdf}
%		\caption{The 2 point torus model of the square lattice.  The fundamental domain is shown at center with copies bordering it.  The two points are shown as the white and black vertices.}
%		\label{Fig:LatticeTorus2ptSq}
%	\end{figure}

	\begin{figure*}[htbp]
		\center
		\includegraphics[width = 1.0\linewidth]{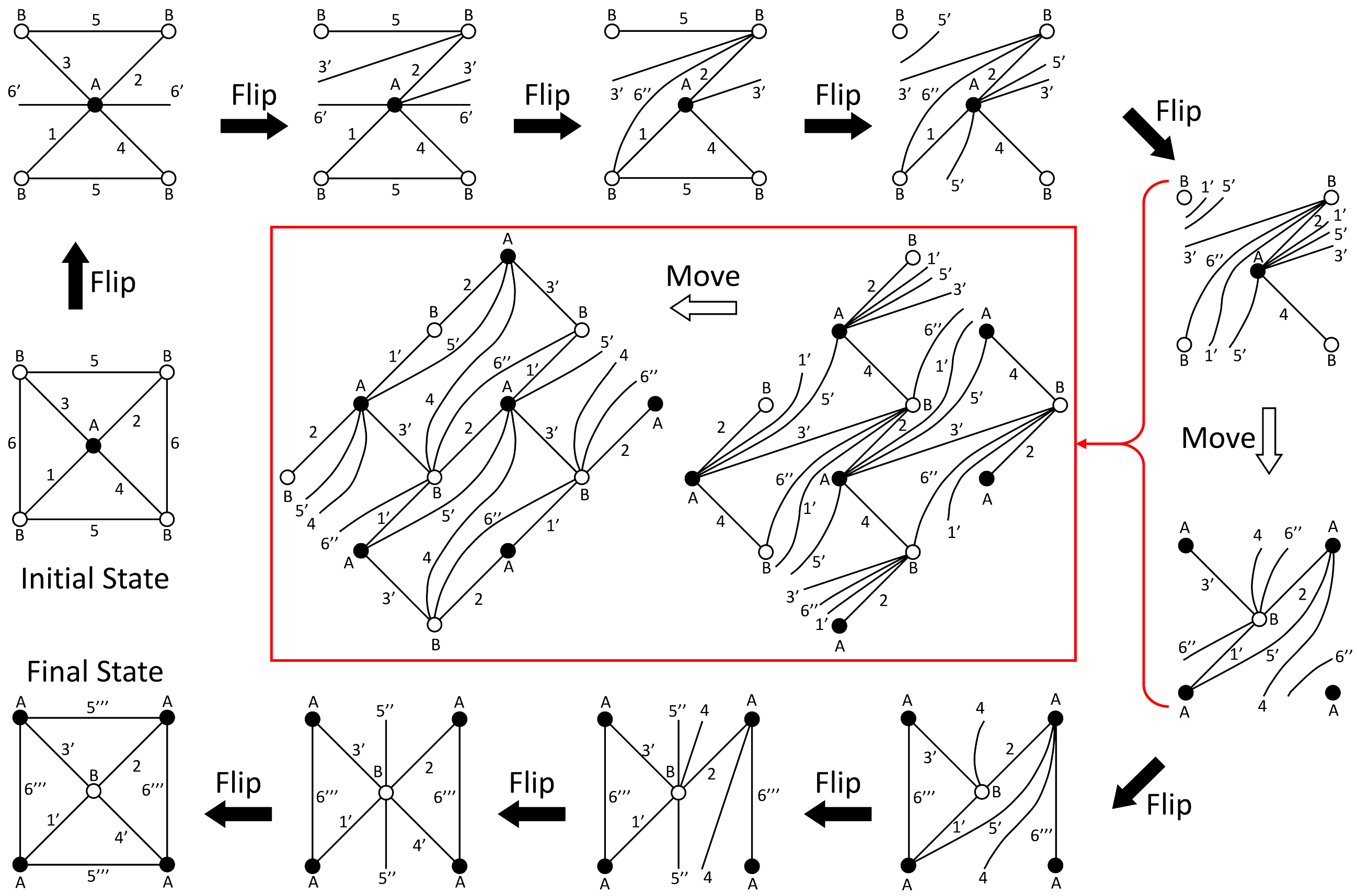}
		\caption{This flip-chart shows the sequence of flips and point movements which take the initial triangulation (center left) to the final triangulation (bottom left) under the action of the braid generator, $\sigma_2$.  The boxed inset shows a larger domain for the point movement step.  To determine the action of the generator on the intersection coordinates, we compare the final and initial triangulation and write down the intervening flips, eq.~\ref{Eq:Flips}. }
		\label{Fig:FlipChart2ptSq}
	\end{figure*}

The procedure from the last section allows us to find the topological entropy per operation for braids defined on graphs which are embedded in compact orientable surfaces.  Here we implement this idea for the two point regular embedding of the square lattice graph on the torus.  As we have seen, the torus embeddings for a lattice graph represent a restriction of possible braid operations to those comprised of a spatially periodic pattern of braid generators.  This particular embedding has the smallest possible fundamental square domain while also allowing for non-trivial braids.  The fundamental square domain, with opposite edges identified, for this torus is shown against the background lattice graph in fig.~\ref{Fig:FundamentalDomains} (upper left).

\begin{figure}[htbp]
		\center
		\includegraphics[width = \linewidth]{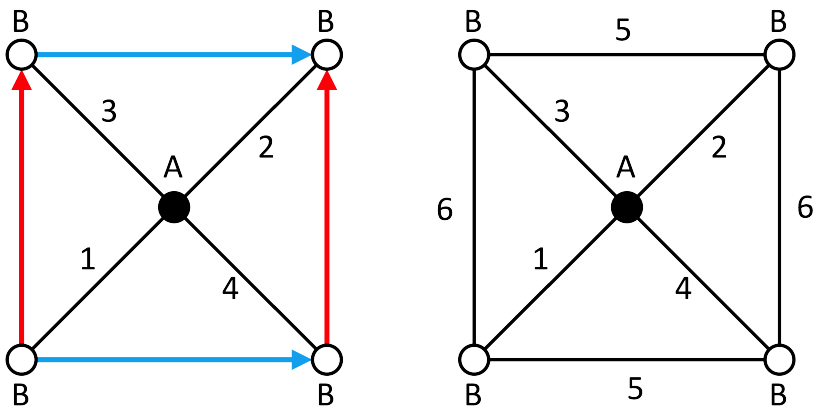}
		\caption{Left: The fundamental domain for the two point torus model of the square lattice.  The blue and red arrows remind us that the top and bottom edges are identified, as are the left and right edges.  The four labeled edges correspond to the braid generators.  Right: The two additional labeled edges make this a triangulation.  These six labeled edges correspond to the intersection coordinates used to encode the train track.}
		\label{Fig:GenAndTriangulation2ptSq}
	\end{figure}

There are four positive braid generators, $\{\sigma_1, \sigma_2,\sigma_3,\sigma_4\}$, (and four inverses, $\sigma^{-1}_i$) corresponding to CCW switches (resp. CW switches) of the two points about the four labeled edges in this graph, (see left side of fig.~\ref{Fig:GenAndTriangulation2ptSq}).  Since there are only two points, only one braid generator at a time can be executed, and there is no distinction between braid generators and braid operations.  We now build up braid words using an alphabet of these 8 braid operations.

In order to set up the intersection coordinate system, we put a triangulation on the torus.  This triangulation is constructed from the four existent graph edges and two additional edges, (see right side of fig.~\ref{Fig:GenAndTriangulation2ptSq}).  Train tracks or closed curves are now uniquely represented by the collection of non-negative integers  $\vec{E} = \left(E_1,E_2,\cdots, E_6\right)$, where $E_i$ counts the minimum number of intersections between the loop and the $i^{th}$ edge in the triangulation.  Note that some indices are shared between the braid generator labeling and the intersection coordinates; the use should be clear from context.

We now want to find the action of our braid generators on the loop intersection coordinates.  We will specify the action of one generator, $\sigma_2$, and give the action of the rest though conjugation with symmetries of the triangulation.  Consider the sequence of flips shown in fig.~\ref{Fig:FlipChart2ptSq}.  As a reminder, the formula for the new edge after a flip is $E' = \max(A+C,B+D) - E \equiv \Delta(A,B,C,D;E)$, where $E$ is the edge between the two triangles, and  $A,B,C,D$ are the edges of the quadrilateral in cyclic order (see fig.~\ref{Fig:TriangleFlip}).  
	
The sequence of flips shown in Fig.~\ref{Fig:FlipChart2ptSq} provide the rules to update an initial train track, given by intersection coordinates $\vec{E} = \left(E_1,E_2,\cdots, E_6\right)$, due to the action of braid generator $\sigma_2$.  The updated coordinates, $\vec{E}^u = \sigma_2 \vec{E}$, are $\vec{E}^u = \left(E'_1,E_2,E'_3,E'_4,E'''_5,E'''_6\right)$, where
\begin{equation}
\begin{split}
E'_6 &= \Delta(E_1,E_3,E_2,E_4;E_6) \\
E'_3 &= \Delta(E_2,E_5,E_2,E'_6;E_3) \\
E'_5 &= \Delta(E_4,E_1,E_2,E'_3;E_5) \\
E''_6 &= \Delta(E_4,E_1,E_2,E'_3;E'_6) \\
E'_1 &= \Delta(E_2,E'_5,E_2,E''_6;E_1) \\
E''_5 &= \Delta(E_2,E'_1,E'_3,E_4;E'_5) \\
E'''_6 &= \Delta(E_2,E'_1,E'_3,E_4;E''_6) \\
E'_4 &= \Delta(E_2,E''_5,E_2,E'''_6;E_4) \\
E'''_5 &= \Delta(E_2,E'_4,E'_1,E'_3;E''_5)
\end{split}
\label{Eq:Flips}
\end{equation}

To get the action of the remaining braid generators, we define three symmetry operators, which act on the triangulation: $R$ for a CCW rotation by $\pi/2$ about the center point (A), $R^{-1}$ for a CW rotation by $\pi/2$, and a mirror inversion $M$ about the ``/" diagonal line through the central point A.  For $i \in [1,2,3,4,5,6]$, we give the action of each symmetry by the permutation $\pi(i)$, where e.g. $[R\vec{E}]_i = E_{\pi(i)}$.   The permutations are defined for the six triangulation edges, but the restriction to the first four edges give the action of the symmetries on the braid generators.
\begin{equation}
\begin{split}
R: \pi &= (3,4,2,1,6,5) \\
R^{-1}: \pi & = (4,3,1,2,6,5) \\
M: \pi &= (1,2,4,3,6,5)
\end{split}
\label{Eq:Symmetries}
\end{equation}
Conjugating $\sigma_2$ with $M$ gives its inverse: $\sigma^{-1}_2 = M\sigma_2 M$ (note that $M^{-1} = M$).  The remaining CCW generators can be readily obtained by conjugation with the rotation operators (the order of operator actions is right to left):
\begin{equation}
\begin{split}
\sigma_1 &= R^{-1}R^{-1} \sigma_2 RR \\
\sigma_3 &= R\sigma_2 R^{-1} \\
\sigma_4 &= R^{-1}\sigma_2 R
\end{split}
\label{Eq:SymConj}
\end{equation}
The remaining CW operators are obtained from eq.~\ref{Eq:SymConj} by replacing $\sigma_2$ with its inverse, $\sigma^{-1}_2$.

We have defined the action of each of the 8 braid operations on the intersection coordinates.  To computationally find the topological entropy per operation (TEPO), we start with an initial set of intersection coordinates, $\vec{E}_0 = (2,2,1,1,4,1)$, which correspond to a closed curve that spirals around the torus four times in one direction for every one time in the other direction.  For pseudo Anosov braids, the specific initial intersection coordinates do not matter, as all non-trivial loops will lead to the same stretching rate.  For a given braid $\beta$, composed of $N$ braid operations, we can use equations~\ref{Eq:Flips}-\ref{Eq:SymConj} to find the updated intersection coordinates, $\vec{E}_1 = \beta \vec{E}_0$, and indeed a sequence of intersection coordinates $(\vec{E}_0,\vec{E}_1, \cdots, \vec{E}_k, \cdots)$.  Define $W_k = \Sigma_i [\vec{E}_k]_i$, as the sum of the intersection coordinates for the $k^{th}$ set of coordinates in this sequence.  Then we can define the TEPO, $ \overline{h}$, for $\beta$ as
\begin{equation}
 \overline{h} = \lim_{k \rightarrow \infty}  \frac{1}{N} \log(\frac{W_k}{W_{k-1}}).
\end{equation}
Computationally, we calculate $ \overline{h}_k = \frac{1}{N}\log(\frac{W_k}{W_{k-1}})$, and stop iterating after $| \overline{h}_k - \overline{h}_{k-1} |$ drops below a desired tolerance.

\section{Search for Maximum TEPO Braids}
\label{Sec:Search}

\begin{table*}[t]
  \centering
  \begin{tabular}{|c|c||c|c|c|c|c|c|}
	 \hline
	 \multicolumn{2}{|c||}{Lattice Types} & \multicolumn{2}{c|}{Square} & \multicolumn{2}{c|}{Triangular} & \multicolumn{2}{c|}{Hexagonal} \\
	 \hline
	 \multicolumn{2}{|c||}{Torus Model} & 2 pt. & 4 pt. & 3 pt. & 4 pt. & 2 pt. & 6 pt. \\
	 %\cline{1-2}
	 %\hhline{|~|~||=|=|=|=|=|=|}
	  \hline
	  \hline
	 \multicolumn{2}{|c||}{\makecell{Number Of \\ Braid Generators}} & 8 & 16 & 18 & 24 & 6 & 18 \\
	 \hline
	 \multicolumn{2}{|c||}{\makecell{Number Of \\ Braid operations}} & 8 & 32 & 18 & 48 & 6 & 48 \\
	 \hline
	 \multirow{3}{*}{\makecell{Max TEPO For \\ Braids Of Length}}
	 & 2 & 0.881373587 & 0.881373587 &  0.881373587 & 0.881373587  & 0.881373587  &  0.881373587 \\
	 \cline{2-8}
	 & 3 & 0.962423650 & 0.962423650  &  0.962423650 & 0.962423650  & 0.962423650 &  0.962423650 \\
	 \cline{2-8}
	 & 4 & 1.061275062 & 1.061275062 &  0.909223230  & 1.061275062  & 0.909223230 &  0.909223230 \\
	 \hline
	  \multirow{2}{*}{\makecell{Checked Braids \\ Up To Length}}
	 %\multicolumn{2}{|c||}{\makecell{Checked Braids \\ Up To Length}} 
	 & A & 9 & 5 & 6 & 4 & 10 & 4 \\
	  \cline{2-8}
	 & B & 10 & 6 & 7 & 6 & 13 & 5 \\
	 \hline
	  \multicolumn{2}{|c||}{\makecell{Max TEPO}} & 1.061275062 & 1.061275062 &  0.962423650  & 1.061275062  & 0.962423650&  0.962423650 \\
	  \hline
	  \multicolumn{8}{|c|}{ $\phi = \frac{1+\sqrt{5}}{2}$, $\delta_s = 1+ \sqrt{2}$, $\log(\delta_s) = 0.881373587$, $2 \log(\phi) = 0.962423650 $, $\log(\phi + \sqrt{\phi}) = 1.061275062$ } \\
	 \hline
	 
  \end{tabular}
  \caption{Max TEPO braid search results.  For the three lattice types and six torus models, we list the following quantities: the number of braid generators (2x the number of graph edges), the number of braid operations (the alphabet that our braid words are composed of), the max TEPO found for braids of short length, the maximum length of braid words checked (A - exhaustively, B - with some common-sense restrictions, see text for discussion), and finally the max TEPO found.  Some useful algebraic representations for TEPO values are also shown.}
  \label{tab:Results}
\end{table*}

For both torus models of each lattice graph, we have an algorithm (see section~\ref{Sec:Example} and Appendix~\ref{Appendix:OtherLatticeGraphs}) which computes the TEPO for a given braid.  We turn now to finding the maximum value of TEPO possible and the lattice graph braid which achieves it.

We proceed by exhaustively testing each braid word of a given length.  If there are $N_o$ operations for a given lattice graph model, then there are $N_o^k$ braid words composed of $k$ operations.  The number of braid words to test becomes untenable as we increase the braid word length, and there is a practical upper limit to this length, set by the acceptable computational time.  Of course, many of the braid words would represent the same braid, but as we do not have the full set of braid relations, we have not sought to reduced the set of braid words this way.  However, we can reduce this set by considering how different braids are connected by the symmetries of the lattice.  Any two braids related by such symmetries will have the same TEPO, and only one will need its TEPO computed.  Instead of mapping out the full symmetry relations between all braid words, we simply restrict the first operation in each braid word to a minimal set of braid operations which generate all braid operations through symmetries.  For example, for the 8 operations of the 2 points square lattice case, see fig.~\ref{Fig:GenAndTriangulation2ptSq}, $\sigma_1$ and $\sigma_1^{-1}$ are related to the rest through rotations.  Furthermore, $\sigma_1$ is related to its inverse via a mirror inversion symmetry along the diagonal.  So, we can start all braid words, in this case, with $\sigma_1$ as the first braid operation.  This reduces the overall run-time by a factor of 8, and allows us to test longer braid words.  The maximum lengths of braid words that we have checked in this manner for each lattice graph model are listed in table~\ref{tab:Results} (``Checked Braids Up To Length", row A).  We further reduced the number of braid words to check by requiring that adjacent braid operations in a braid word not share any edges in the graph.  This has the effect of excluding braid words that contain adjacent braid generators which are each other's inverses, as well as excluding repeated braid generators.  While this reduced set does exclude some unique braid words (these have non-maximal TEPO), it allows us to check longer braid words, see (``Checked Braids Up To Length", row B, in table~\ref{tab:Results}). 

Table~\ref{tab:Results} also contains the max TEPO values found for each lattice graph model. There is a clear maximum TEPO value of $\overline{h} = 1.061275062$, which appears for square lattices and the 4 point triangular lattice, starting with braids of length 4.  This persists as the maximum TEPO value through all braid words that we have checked.  In the next section, we will have a close look at the braids which realize this max TEPO value.

\section{Candidate Maximum TEPO Braid}
\label{Sec:MaxTEPOBraid}

\begin{figure}[htbp]
		\center
		\includegraphics[width = \linewidth]{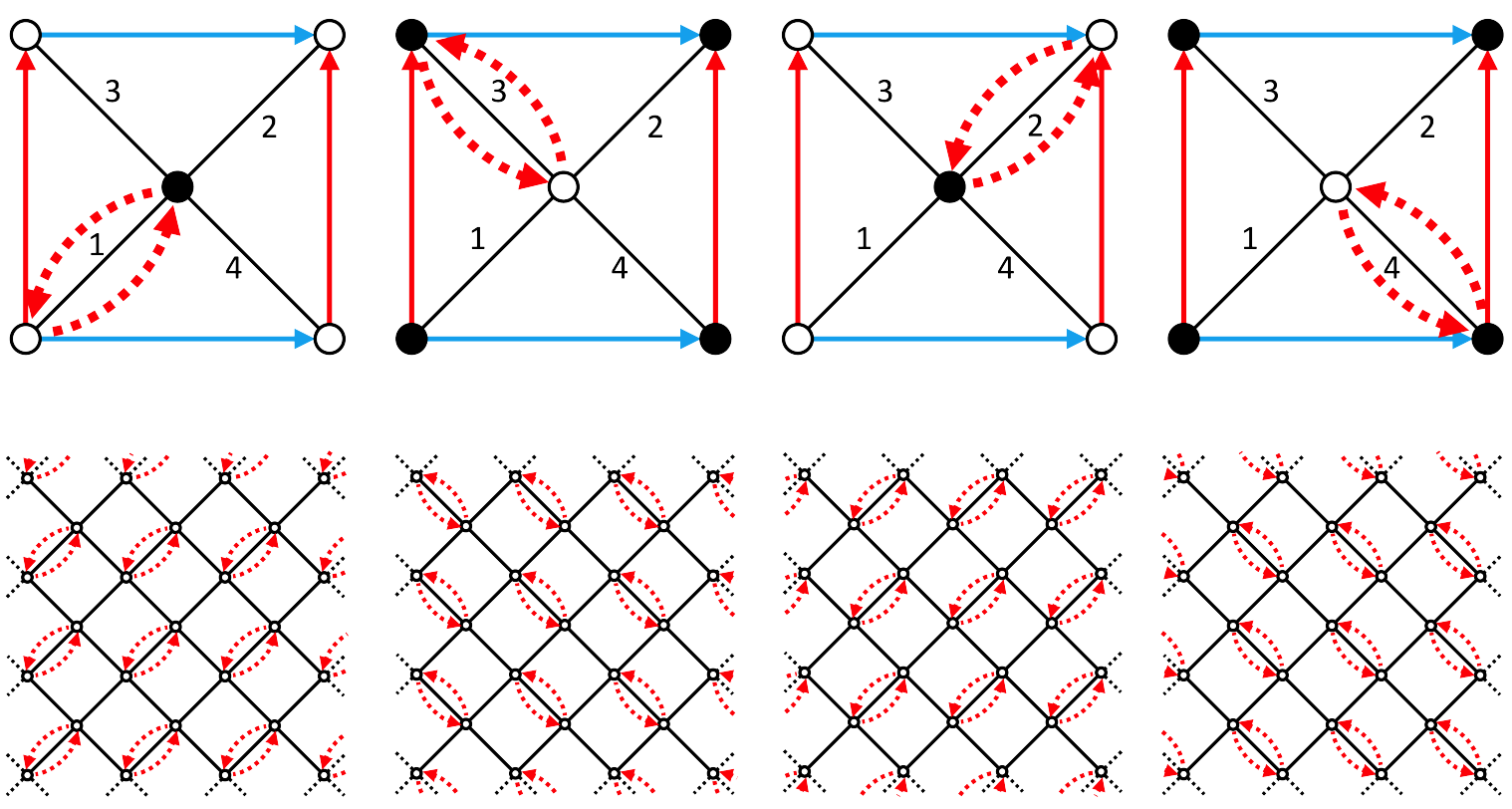}
		\caption{The max TEPO braid, eq.~\ref{eq:maxTEPObraid}, shown as point-pair switches (left to right).  The top row shows this on the two point torus model of the square lattice graph, while the bottom row depicts this on a larger domain.
		}
		\label{Fig:maxTEPObraid}
	\end{figure}

We searched through hundreds of millions of braid words (see table~\ref{tab:Results}) and found a very small number of braid words that attained the maximum TEPO value of $\overline{h} = 1.061275062$.  Furthermore, all such braid words on a given lattice graph model are related by symmetries (rotations and time-reversals), and can be considered the same braid.  Though these braids are different for the 2 point square, 4 point square, and 4 point triangular lattice graph models (e.g. incompatible number of moving points), there is a way in which they can all be considered essentially the same braid.  If you extend each of these lattice graph model braids to the full lattice in $\mathbb{R}^2$ through the periodicity of the torus (e.g. see the bottom of fig.~\ref{Fig:maxTEPObraid}), then the braids which result are isotopic.  Let us refer to this general max TEPO braid as $\beta^{*}$.  In this section we will investigate a particular braid word example of $\beta^{*}$ (see eq.~\ref{eq:maxTEPObraid}).  However, the properties that result are either universal to $\beta^{*}$ or are easy to generalize to other braid word instantiations of $\beta^{*}$.  

In terms of the braid generators introduced in section~\ref{Sec:Example} for the 2 point square lattice model, this braid can be written as
\begin{equation}
    \beta^{*}_{2PtSq} = \sigma_1 \sigma_3 \sigma_2 \sigma_4.
    \label{eq:maxTEPObraid}
\end{equation}
For completeness, we also give braid words corresponding to $\beta^{*}$ on the 4 point square and triangular lattice models:
\begin{equation}
\begin{split}
    \beta^{*}_{4PtSq} &= (\sigma_1 \sigma_9)(\sigma_4 \sigma_6)(\sigma_3 \sigma_8)(\sigma_2 \sigma_7) \\
    \beta^{*}_{4PtTri} &= (\sigma_1 \sigma_4)(\sigma_3 \sigma_9)(\sigma_2 \sigma_{12})(\sigma_4 \sigma_{10}).
    \label{eq:maxTEPObraid2}
\end{split}
\end{equation}
Here, pairs of generators comprising a braid operator are grouped with parentheses and the indexing for braid generators is given by the labeled edges in the relevant diagrams of fig.~\ref{Fig:AllCasesTriangulations} in appendix~\ref{Appendix:OtherLatticeGraphs}.

We now focus on $\beta^{*}_{2PtSq}$.  The only other braid words that give the same TEPO value on this lattice graph model are those related by cyclic shifts (e.g. $\sigma_2 \sigma_4 \sigma_1 \sigma_3$) or time reversal symmetry (e.g. $\sigma_4^{-1} \sigma_2^{-1} \sigma_3^{-1} \sigma_1^{-1}$), and therefore constitute essentially the same braid.  This braid is rotational symmetric, and rotating the underlying lattice by $\pi/2$ or separately by $-\pi/2$, see eq.~\ref{Eq:Symmetries}, is equivalent to a cyclic shift of the braid word.  Indeed, if $R$ represents a CCW rotation of the lattice by $\pi/2$, then we can represent this braid by alternating rotations with the braid generator $\sigma_1$: $\beta^{*}_{2PtSq} = \sigma_1 R \sigma_1 R \sigma_1 R \sigma_1 R$ (read left to right).

The series of pair-wise switches which make up this braid can be seen in fig.~\ref{Fig:maxTEPObraid}.  On the top, the max TEPO braid is realized on the 2 point torus model of the square lattice.  On the bottom we show this braid over a larger portion of the square lattice.  Here the spatial periodicity is more visually manifest, and we can see that this braid is relatively simple.  

\begin{figure}[htbp]
		\center
		\includegraphics[width = \linewidth]{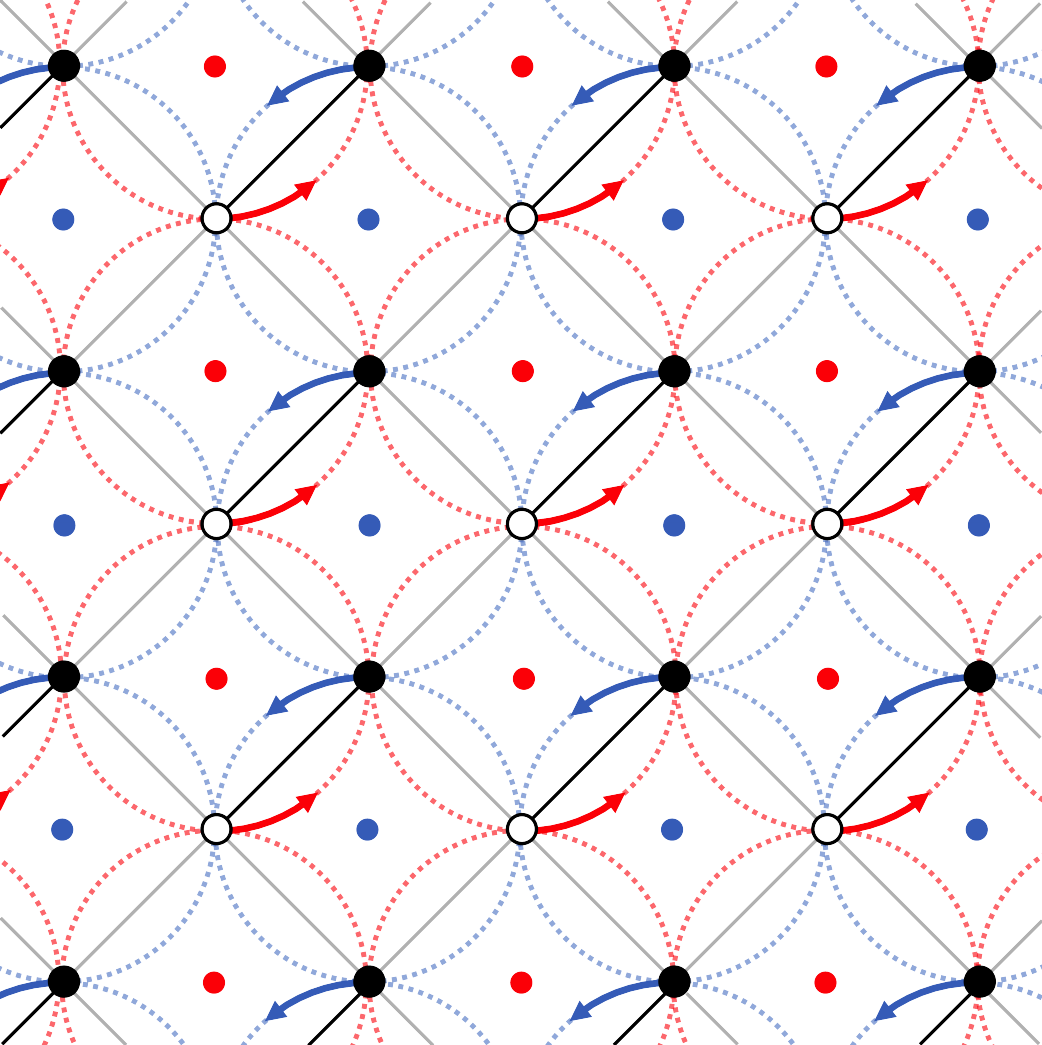}
		\caption{One possibility for geometric motion of points compatible with the max TEPO braid.  The red and blue dots mark the centers for circular motion of the white and black points, respectively.  The doted lines mark the paths and the arrows show the direction of motion.  Note that the red and blue sets of circular oscillators are a half rotation out of phase from one another.}
		\label{Fig:GeometricMovement}
	\end{figure}

This representation of the braid is still purely algebraic, and it is helpful to articulate a simple model of specific geometric motions that have this braiding topology.  First note that if you follow the path of any one point, it is periodic after one application of the braid, making this a pure braid (braids for which the induced permutation of strands is the identity).  Furthermore this point traverses the sides of one of the square lattice faces (CCW).  We simply make this circular motion.  Now we have a set of circular oscillators, each centered on the faces of the square lattice, all of which are moving CCW.  These oscillators are separated into two groups by a checker-board coloring of this square lattice, and the rotation of each group is out of phase by $\pi$ with one another.  See fig.~\ref{Fig:GeometricMovement} for a view of this movement.

\begin{figure}[htbp]
		\center
		\includegraphics[width = \linewidth]{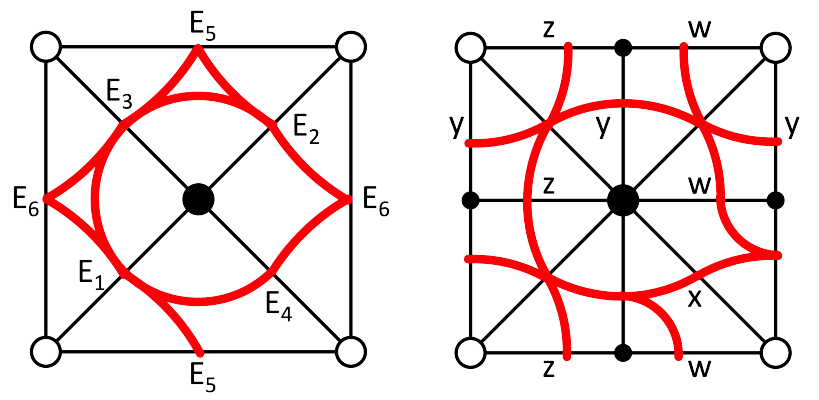}
		\caption{Left: the invariant measured train track for $\beta^*$ (in red).  Right: the same invariant measured train track, but represented on a new triangulation - required for the Veering triangulation splitting sequence.}
		\label{Fig:InvariantTrainTracks}
	\end{figure}

Repeated action of the maximum TEPO braid on any initial, non-trivial train track, along with normalizing the measure (dividing by the largest intersection coordinate), will lead to a sequence of train tracks which asymptotically approach the projective invariant measured train track for $\beta^*_{2PtSq}$.  This invariant measured train track can be seen on the left side of fig.~\ref{Fig:InvariantTrainTracks}, and has normalized intersection coordinates $\vec{E} = $(1.0, 0.57230276, 0.89005364, 0.34601434, 0.65398566, 0.91831709).  The action of $\beta^*_{2PtSq}$ on these coordinates give
\begin{equation}
    \beta^*_{2PtSq} \vec{E} = \lambda \vec{E},
    \label{eq:tteigen}
\end{equation}
where the braid dilation, $\lambda = \exp(4\overline{h}) = 69.7627534$, is obtained from the TEPO value, $\overline{h} = 1.061275062$.  Eq.~\ref{eq:tteigen} is an eigenvalue equation for $\beta^*_{2PtSq}$, considered as an operator.  Indeed, given the intersection coordinates of the invariant measured train track, i.e. the eigenvector $\vec{E}$, we can replace the max operation in each triangulation flip, eq.~\ref{Eq:Flip}, with the corresponding linear operator.  All the flips in each braid generator of our max TEPO braid then combine to give a single matrix with $\lambda$ as its largest eigenvalue.  This method can be used to find an analytical form for the TEPO value.  However, we will use another method to find the analytical form for $\overline{h}$, which will additionally provide us with a total conjugacy invariant of the braid $\beta^*_{2PtSq}$.

We will construct a Veering triangulation~\cite{Agol2011IdealTO} for the max TEPO braid.  The full definition of a Veering triangulation concerns 3D triangulations of the mapping torus of closed surface mapping classes, but we will not need the full machinery.  For our purposes, a Veering triangulation consists of a periodic sequence of 2D triangulations of our torus, related by triangulation flips, which admit a certain edge coloring structure.  We should be able to color each edge black or blue such that for each edge flip, the four non-diagonal edges of the associated quadrilateral alternate colors when ordered cyclically, and the two black edges form a ``Z" with the pre-flipped diagonal edge.  

	\begin{figure*}[htbp]
		\center
		\includegraphics[width = 1.0\linewidth]{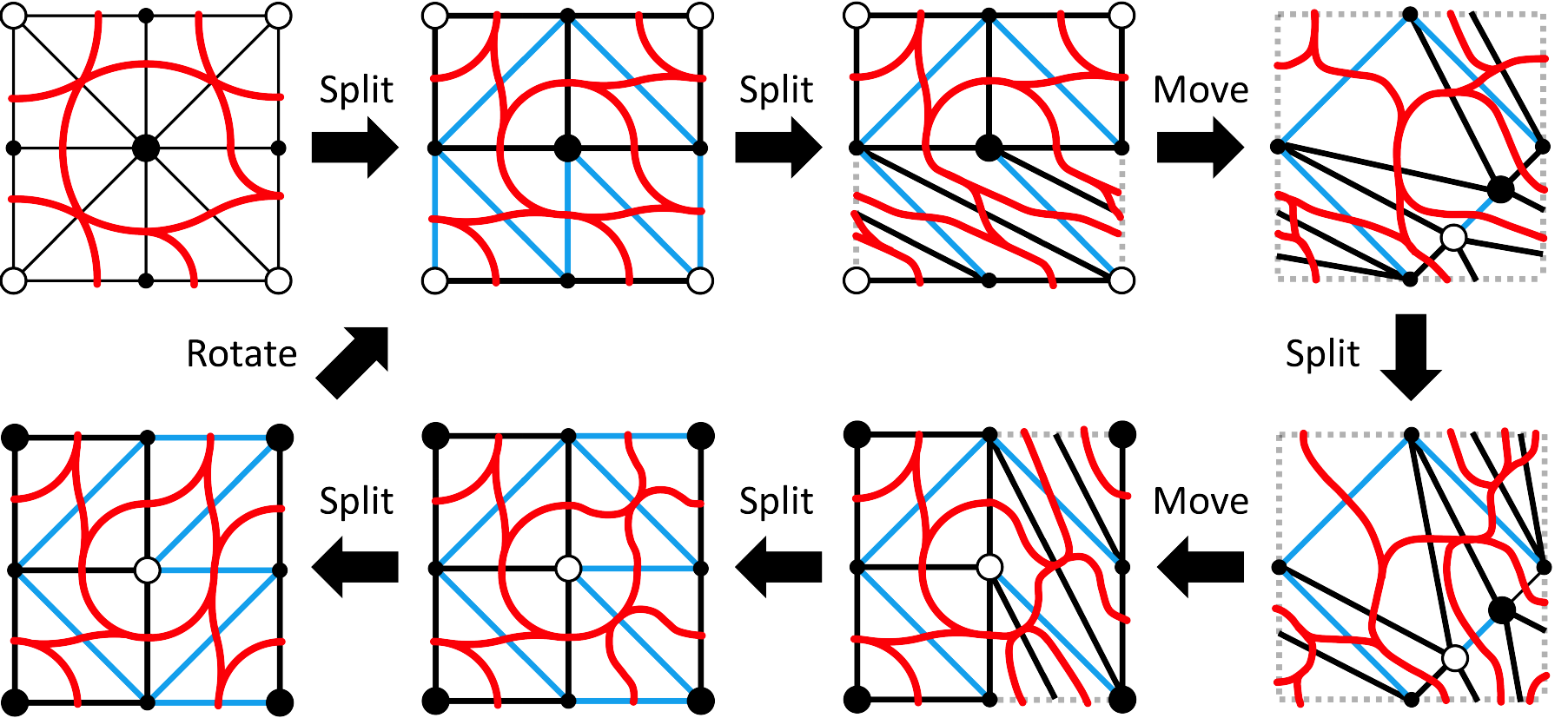}
		\caption{The eventually periodic train track splitting sequence for $\beta^*_{2PtSq}$.  The Veering triangulation structure is given by the black/blue coloring of edges.}
		\label{Fig:SplittingSequence}
	\end{figure*}

Fortunately, there is an automatic method for generating this structure.  We first must start with a slightly modified invariant train track.  This procedure requires that each component of the surface, cut by the invariant train track, contain a vertex of the triangulation.  We have two one-cusp regions which already have a vertex (the black and white circles in the left part of fig.~\ref{Fig:InvariantTrainTracks}), and two three-cusp regions without vertices.  We introduce two new vertices in these three-cusp regions to get the new triangulation and train track seen on the right side of fig.~\ref{Fig:InvariantTrainTracks}.  Here the train track measure is completely specified by just four intersection coordinates, $(w,x,y,z) = (0.16811739, 0.34601434,  0.40418537, 0.48586828)$.  The other intersection coordinates can be found from the train track switch conditions.  Note that each triangle in this triangulation now satisfies a triangle equality, and therefore contains only one switch (train track vertex).

Starting with this new triangulation and its intersection coordinates, we generate the next triangulation in the sequence by flipping the edge with the largest intersection coordinate (whose associated train tracks operation is a split).  Repeating this procedure, we obtain a splitting sequence of triangulations, as seen in fig.~\ref{Fig:SplittingSequence}, which is eventually periodic.  As our braid has rotational symmetry, we can stop the splitting sequence after $1/4^{th}$  of its length and connect beginning and ending triangulations by a rotation.  This symmetry reduction will simplify the analysis.  Note that we have also colored the edges of each triangulation in the eventually periodic portion of the splitting sequence with black and blue in the manner described above.  Thus, this sequence of triangulations is a total conjugacy invariant, and could be used to definatively differentiate this pseudo-Anosov braid from others which might have the same topological entropy (though we have found no others in our survey).

	\begin{figure}[htbp]
		\center
		\includegraphics[width = \linewidth]{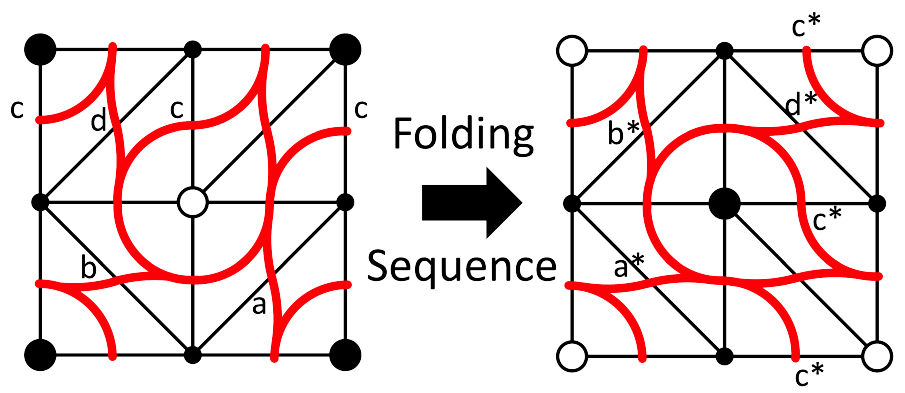}
		\caption{The beginning and end states of the folding sequence.  The Perron-Frobenius transition matrix, eq.~\ref{eq:transitionMat}, relates the final coordinates (starred) to the initial coordinates.}
		\label{Fig:FoldingSequence}
	\end{figure}

If we reverse the splitting sequence we have a folding sequence, fig.~\ref{Fig:FoldingSequence}, from which we can construct the Perron-Frobenius matrix, $\mathbf{A}$, used to update the intersection coordinates.  In particular, we construct this transition matrix for a minimal generating subset of intersection coordinates, $\vec{E}_m = (a,b,c,d)$ on  the left side of fig.~\ref{Fig:FoldingSequence} (all other intersection coordinates can be calculated from this minimal set, and the train track switch conditions).  Each flip (corresponding to a train track folding move) in the sequence gives a linear operator on the intersection coordinates.  Accumulating the effect of each flip and the permutation needed to map the final train track onto the initial one, gives the new coordinates, $\vec{E}^*_m  = \mathbf{A} \vec{E}_m$, where
\begin{equation}
    \mathbf{A} = 
    \begin{bmatrix}
    0 & 1 & 0 & 0 \\
    0 & 0 & 0 & 1 \\
    0 & 1 & 1 & 1 \\
    1 & 1 & 2 & 1 
    \end{bmatrix}
    \label{eq:transitionMat}
\end{equation}	
This transition matrix has a characteristic polynomial of $\lambda^4 -2\lambda^3 - 2\lambda^2 - 2\lambda + 1$.  The largest root of this polynomial is
\begin{equation}
    \lambda = \phi + \sqrt{\phi},
    \label{eq:analyticaldilation}
\end{equation}
where $\phi = (1+\sqrt{5})/2$ is the golden ratio.  The TEPO value is therefore $\overline{h} = \log(\phi + \sqrt{\phi})$.  Curiously, this characteristic polynomial is palindromic, and $\phi + \sqrt{\phi}$ is a Salem number, algebraic numbers with particularly interesting properties.  

In a satisfying connection, $\phi + \sqrt{\phi}$ is directly related to the braid dilation of a particular historical taffy stretching device (see table 1 from reference~\cite{TaffyThiffeault}, and fig. 16 from the arxiv version of that article).  Indeed, a potential mixing device (which we introduce in fig.~\ref{Fig:Mixer}) based on the max TEPO braid, has point movements that are the natural 2D generalization of the rod movement in this taffy stretching device.

\begin{figure}[htbp]
		\center
		\includegraphics[width = \linewidth]{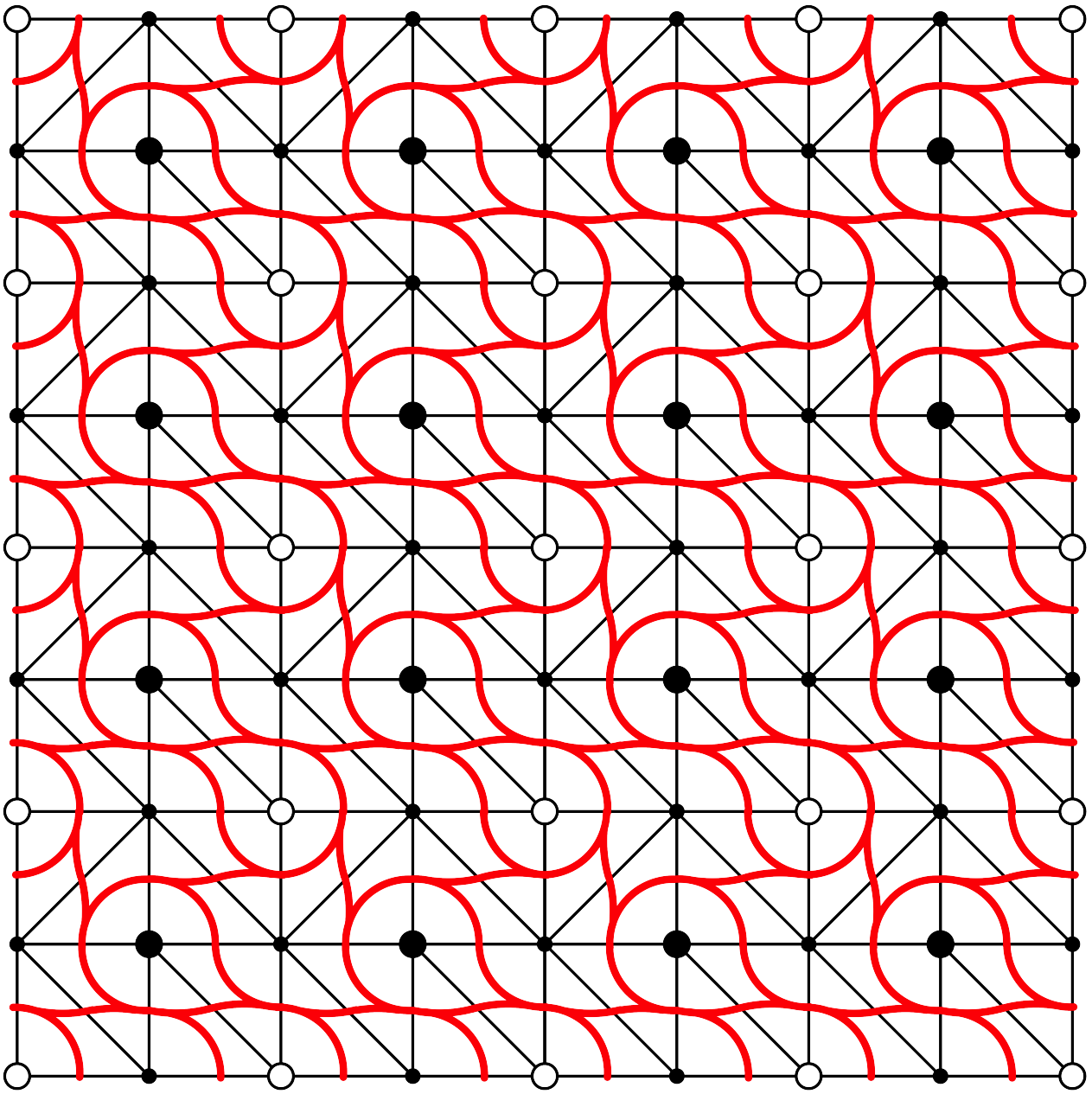}
		\caption{And extended section of the invariant train track for $\beta^*$, in red.  Note that some edges in the original square lattice graph are not present, as this triangulation arises from the Veering structure.}
		\label{Fig:BigInvTrainTracks}
	\end{figure}

All of these properties, the Veering triangulation, the invariant projective measured train track, the transition matrix, and the TEPO value, naturally extend to the full square lattice, just as $\beta^*_{2PtSq}$ can be thought of as the braid $\beta^*$ on this larger domain, see fig.~\ref{Fig:maxTEPObraid}.  Fig.~\ref{Fig:BigInvTrainTracks} shows the invariant train track for a larger section of the square lattice.  This braid can also be found on triangular lattices, though confined to a square lattice sub-graph.  

Given the computational evidence we have accumulated and the relative simplicity of the braid $\beta^*$, we conjecture that $\beta^*$ is the unique braid which maximizes the topological entropy per operation on all planar lattice graphs:

\begin{conj}
Let $\mathbf{G}$ be the countably infinite set of regular maps on the torus (i.e. $\mathcal{G} \in \mathbf{G}$ is an embedded graph in the torus that is arc-transitive), let $\mathbf{B}_{\mathcal{G}}$ be the set of braid words formed from any number of braid generators defined on $\mathcal{G}$, and let $\overline{h}(\beta)$ be the topological entropy per operation of the braid word $\beta \in \mathbf{B}_{\mathcal{G}}$.  Then
\begin{equation}
    \max_{\beta \in \mathbf{B}_{\mathcal{G}}, \mathcal{G} \in \mathbf{G}} \overline{h}(\beta) = \log(\phi+\sqrt{\phi}),
\end{equation}
where $\phi =  \frac{1+\sqrt{5}}{2}$.  Furthermore, when extended to the plane using the periodicity of the pertinent torus model, all braid words $\beta$, for which $\overline{h}(\beta) = \log(\phi+\sqrt{\phi})$, represent the same braid up to conjugacy and time-reversal symmetry.
\end{conj}

\section{Connections and Conclusions}
\label{Sec:Connections}

We started this paper by highlighting a number of applications that connect the topological entropy of braids with the mixing characteristics of dynamical systems.  Motivated by this, we have put forward a general method for finding the topological entropy and TEPO for graph generated surface braids (of which regular braids are a special case), and we have found numerical evidence for a unique braid which, we conjecture, maximizes the TEPO for all lattice graph braids on the plane.  Circling back, here we consider this particular max TEPO braid in terms of potential applications.

\begin{figure}[htbp]
		\center
		\includegraphics[width = \linewidth]{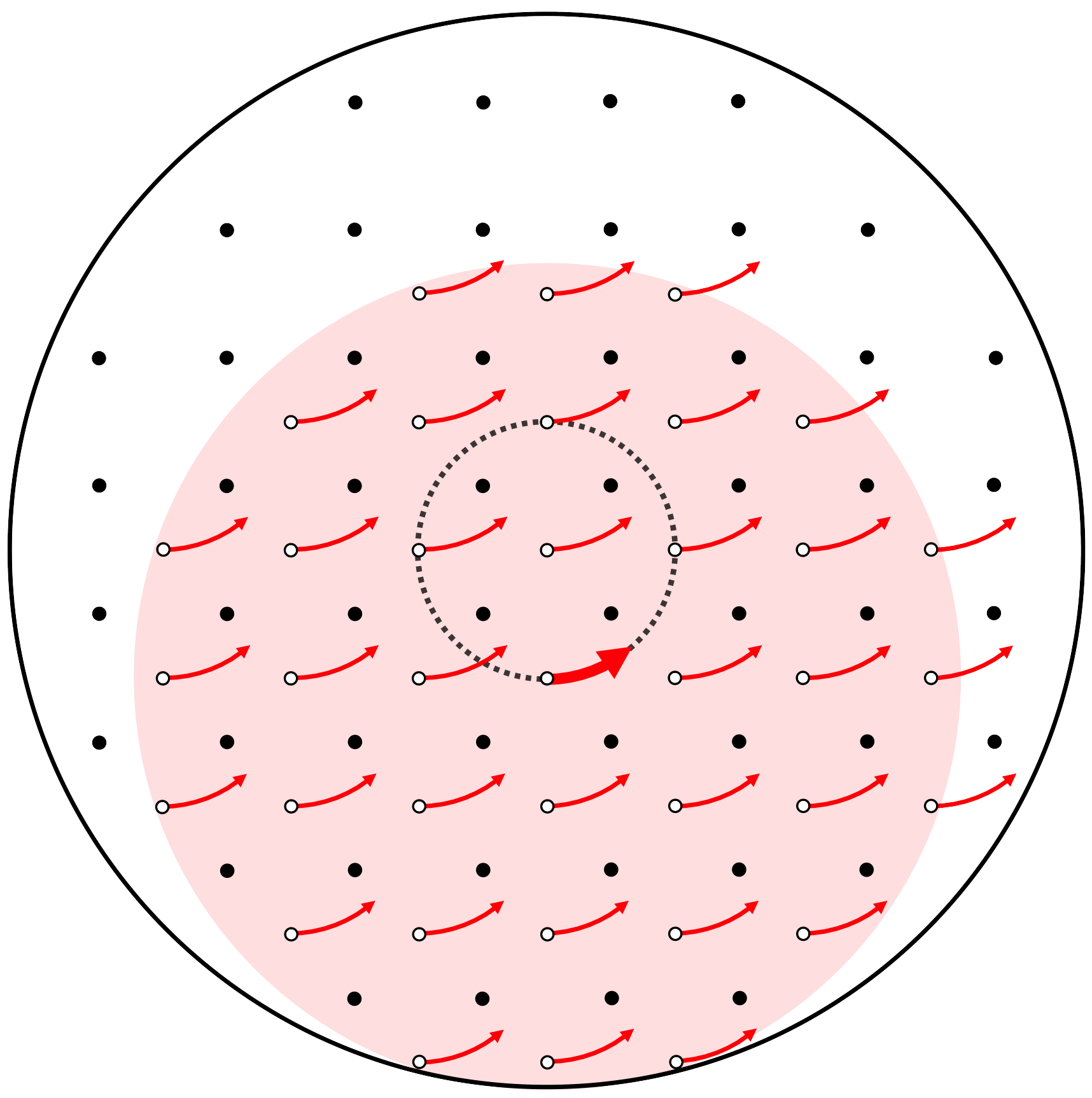}
		\caption{A potential mixing device based on the max TEPO braid.  In this top-down view, the large black circle is the boundary of the cylindrical vessel.  The square array of black dots are mixing rods fixed to the bottom, while the smaller square array of white dots (inside the highlighted disk) are the moving mixing rods affixed to a motor mechanism on the top of the vessel.  The arrows show the rigid irrotational motion of the rods, while the small circle shows the path of one rod (center of the moving array, with bold arrow) over the course of one cycle of the motion.  It is straight-forward to adapt this geometry to more dense arrays of rods.}
		\label{Fig:Mixer}
	\end{figure}

First, the max TEPO braid can be realized as a rod-based stirring device with a particularly simple mechanism that should be easy to implement.  Consider the geometric motion shown in fig.~\ref{Fig:GeometricMovement}. To be clear, we are considering some finite subset of the infinite lattice of points, say the points within a circle of a given radius. These points, now thought of as stirring rods, come in two groups, each arranged as a square lattice (the chess-board coloring subsets of the original square lattice).  Each group moves rigidly, and irrotationally, in a circle of radius half that of the group lattice spacing.  The circular motion of the two groups are half a rotation out of phase with each other.  To realize this rod motion in a cylindrical mixing vessel, one group could be attached at the bottom and the other at the top, each driven by a separate motor mechanism.  However, a simpler realization takes the bottom-anchored mixing rod group as stationary, and drives the top group to execute circles that are twice as large as before, see fig.~\ref{Fig:Mixer}.  This would require only one driving motor, and the same amount of work would be done on the fluid (half the rods moving twice the distance) and therefore the same energy expenditure.  However, assuming similar rod speeds, this would take twice as long for the same mixing effect.

As a second application, we consider the connection between max TEPO braids and the movement of topological defects in systems of 2D active nematic microtubule bundles.  These systems exhibit uniform extensile dynamics, as microtubule bundles slide past each other, driven by micro-molecular motors~\cite{DogicANMT}.  The injection of energy at this small scale leads to a uniform stretching rate for microtubule bundle material curves.  As opposed to rod-based stirring devices, where the braid is specified and the stretching is what results, here the stretching rate is specified and the braids formed from the motion of $+\frac{1}{2}$ topological defects are the emergent feature.  These braids must be pseudo-Anosov and have positive topological entropy.  Indeed, when confined to an annular channel~\cite{C6SM02310J}, the positive defects move in such a way as to generate the silver braid, which maximizes TEPO for annular braids~\cite{MR2317905,MR2861264}.  We have found analogous results for the braiding of four positive defects on a sphere, and maximizing TEPO seems to be a poorly understood tendency of these systems (with the caveat that we have considered configurations where defect pair creation and annihilation events are rare).   This brings up the possibility for our max TEPO lattice braid to be realized by an active nematic microtubule system.  

	\begin{figure}[htbp]
		\center
		\includegraphics[width = \linewidth]{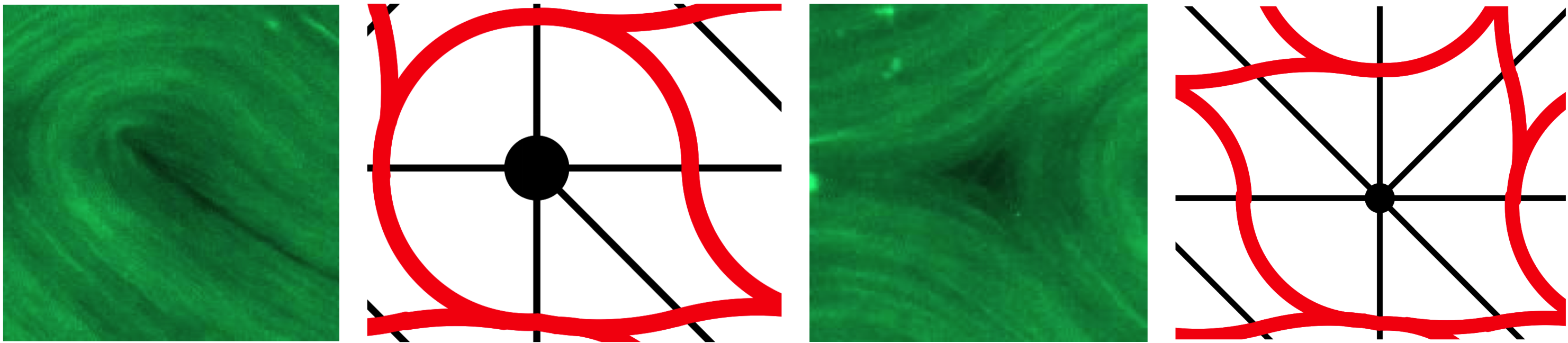}
		\caption{From left to right: an experimental image of a $+\frac{1}{2}$ topological defect, a section of a train-track surrounding a one-prong singularity in the map, an image of a $-\frac{1}{2}$ defect, train-track surrounding a three-prong singularity.  The train-track images are details from fig.~\ref{Fig:BigInvTrainTracks}, while the microtubule bundle images, from the topological chaos in active nematics study~\cite{ActiveNemTopChaos2019}, are courtesy of Amanda Tan and the Hirst lab at UC Merced.}
		\label{Fig:DefectProngCompare}
	\end{figure}

There is a general connection between $+\frac{1}{2}$ defects and the one-cusp regions of the surface cut by train-tracks (also thought of as a one-prong singularity in the corresponding map), as well as between the $-\frac{1}{2}$ defects and the three-cusp regions, see fig.~\ref{Fig:DefectProngCompare}.  In both cases, the train-track structure encodes the topology of how the micotubule bundles bend around the topological defects.  This association is also dynamically meaningful, as it is exclusively the movement of the one-cusp points and the $+\frac{1}{2}$  topological defects that account for the production of topological entropy~\cite{ActiveNemTopChaos2019}.  To realize our max TEPO braid in an active nematic microtubule system, we would likely need to pin down the $-\frac{1}{2}$ defects at lattice sites.  Efforts are underway to find evidence for the existence of the max TEPO braid in an active nematic microtubule system.

%things to include: look at graphs on the sphere (regular maps -  platonic solids) and on higher genus surfaces (hyperbolic surfaces - neg. Gaussian curvature), also adapted the algorithm to deal with arbitrary point motions on the plane.  

%\normalsize
%\textbf{References}\small
\bibliography{MaxTEPOLatticeBraidBib}

\begin{thebibliography}{55}%
\makeatletter
\providecommand \@ifxundefined [1]{%
 \@ifx{#1\undefined}
}%
\providecommand \@ifnum [1]{%
 \ifnum #1\expandafter \@firstoftwo
 \else \expandafter \@secondoftwo
 \fi
}%
\providecommand \@ifx [1]{%
 \ifx #1\expandafter \@firstoftwo
 \else \expandafter \@secondoftwo
 \fi
}%
\providecommand \natexlab [1]{#1}%
\providecommand \enquote  [1]{``#1''}%
\providecommand \bibnamefont  [1]{#1}%
\providecommand \bibfnamefont [1]{#1}%
\providecommand \citenamefont [1]{#1}%
\providecommand \href@noop [0]{\@secondoftwo}%
\providecommand \href [0]{\begingroup \@sanitize@url \@href}%
\providecommand \@href[1]{\@@startlink{#1}\@@href}%
\providecommand \@@href[1]{\endgroup#1\@@endlink}%
\providecommand \@sanitize@url [0]{\catcode `\\12\catcode `\$12\catcode
  `\&12\catcode `\#12\catcode `\^12\catcode `\_12\catcode `\%12\relax}%
\providecommand \@@startlink[1]{}%
\providecommand \@@endlink[0]{}%
\providecommand \url  [0]{\begingroup\@sanitize@url \@url }%
\providecommand \@url [1]{\endgroup\@href {#1}{\urlprefix }}%
\providecommand \urlprefix  [0]{URL }%
\providecommand \Eprint [0]{\href }%
\providecommand \doibase [0]{http://dx.doi.org/}%
\providecommand \selectlanguage [0]{\@gobble}%
\providecommand \bibinfo  [0]{\@secondoftwo}%
\providecommand \bibfield  [0]{\@secondoftwo}%
\providecommand \translation [1]{[#1]}%
\providecommand \BibitemOpen [0]{}%
\providecommand \bibitemStop [0]{}%
\providecommand \bibitemNoStop [0]{.\EOS\space}%
\providecommand \EOS [0]{\spacefactor3000\relax}%
\providecommand \BibitemShut  [1]{\csname bibitem#1\endcsname}%
\let\auto@bib@innerbib\@empty
%</preamble>
\bibitem [{\citenamefont {Boyland}\ \emph {et~al.}(2000)\citenamefont
  {Boyland}, \citenamefont {Aref},\ and\ \citenamefont {Stremler}}]{MR1742169}%
  \BibitemOpen
  \bibfield  {author} {\bibinfo {author} {\bibfnamefont {P.~L.}\ \bibnamefont
  {Boyland}}, \bibinfo {author} {\bibfnamefont {H.}~\bibnamefont {Aref}}, \
  and\ \bibinfo {author} {\bibfnamefont {M.~A.}\ \bibnamefont {Stremler}},\
  }\href {\doibase 10.1017/S0022112099007107} {\bibfield  {journal} {\bibinfo
  {journal} {J. Fluid Mech.}\ }\textbf {\bibinfo {volume} {403}},\ \bibinfo
  {pages} {277} (\bibinfo {year} {2000})}\BibitemShut {NoStop}%
\bibitem [{\citenamefont {Thiffeault}\ \emph {et~al.}(2008)\citenamefont
  {Thiffeault}, \citenamefont {Finn}, \citenamefont {Gouillart},\ and\
  \citenamefont {Hall}}]{MR2464304}%
  \BibitemOpen
  \bibfield  {author} {\bibinfo {author} {\bibfnamefont {J.-L.}\ \bibnamefont
  {Thiffeault}}, \bibinfo {author} {\bibfnamefont {M.~D.}\ \bibnamefont
  {Finn}}, \bibinfo {author} {\bibfnamefont {E.}~\bibnamefont {Gouillart}}, \
  and\ \bibinfo {author} {\bibfnamefont {T.}~\bibnamefont {Hall}},\ }\href
  {\doibase 10.1063/1.2973815} {\bibfield  {journal} {\bibinfo  {journal}
  {Chaos}\ }\textbf {\bibinfo {volume} {18}},\ \bibinfo {pages} {033123, 8}
  (\bibinfo {year} {2008})}\BibitemShut {NoStop}%
\bibitem [{\citenamefont {Thiffeault}\ and\ \citenamefont
  {Finn}(2006)}]{MR2317905}%
  \BibitemOpen
  \bibfield  {author} {\bibinfo {author} {\bibfnamefont {J.-L.}\ \bibnamefont
  {Thiffeault}}\ and\ \bibinfo {author} {\bibfnamefont {M.~D.}\ \bibnamefont
  {Finn}},\ }\href {\doibase 10.1098/rsta.2006.1899} {\bibfield  {journal}
  {\bibinfo  {journal} {Philos. Trans. R. Soc. Lond. Ser. A Math. Phys. Eng.
  Sci.}\ }\textbf {\bibinfo {volume} {364}},\ \bibinfo {pages} {3251} (\bibinfo
  {year} {2006})}\BibitemShut {NoStop}%
\bibitem [{\citenamefont {Thiffeault}(2010)}]{BraidsEntPartTraj}%
  \BibitemOpen
  \bibfield  {author} {\bibinfo {author} {\bibfnamefont {J.-L.}\ \bibnamefont
  {Thiffeault}},\ }\href {\doibase 10.1063/1.3262494} {\bibfield  {journal}
  {\bibinfo  {journal} {Chaos: An Interdisciplinary Journal of Nonlinear
  Science}\ }\textbf {\bibinfo {volume} {20}},\ \bibinfo {pages} {017516}
  (\bibinfo {year} {2010})}\BibitemShut {NoStop}%
\bibitem [{\citenamefont {Gouillart}\ \emph {et~al.}(2006)\citenamefont
  {Gouillart}, \citenamefont {Thiffeault},\ and\ \citenamefont
  {Finn}}]{PhysRevE.73.036311}%
  \BibitemOpen
  \bibfield  {author} {\bibinfo {author} {\bibfnamefont {E.}~\bibnamefont
  {Gouillart}}, \bibinfo {author} {\bibfnamefont {J.-L.}\ \bibnamefont
  {Thiffeault}}, \ and\ \bibinfo {author} {\bibfnamefont {M.~D.}\ \bibnamefont
  {Finn}},\ }\href {\doibase 10.1103/PhysRevE.73.036311} {\bibfield  {journal}
  {\bibinfo  {journal} {Phys. Rev. E}\ }\textbf {\bibinfo {volume} {73}},\
  \bibinfo {pages} {036311} (\bibinfo {year} {2006})}\BibitemShut {NoStop}%
\bibitem [{\citenamefont {Allshouse}\ and\ \citenamefont
  {Thiffeault}(2012)}]{ALLSHOUSE201295}%
  \BibitemOpen
  \bibfield  {author} {\bibinfo {author} {\bibfnamefont {M.~R.}\ \bibnamefont
  {Allshouse}}\ and\ \bibinfo {author} {\bibfnamefont {J.-L.}\ \bibnamefont
  {Thiffeault}},\ }\href {\doibase https://doi.org/10.1016/j.physd.2011.10.002}
  {\bibfield  {journal} {\bibinfo  {journal} {Physica D: Nonlinear Phenomena}\
  }\textbf {\bibinfo {volume} {241}},\ \bibinfo {pages} {95} (\bibinfo {year}
  {2012})}\BibitemShut {NoStop}%
\bibitem [{\citenamefont {Filippi}\ \emph {et~al.}(2020)\citenamefont
  {Filippi}, \citenamefont {Budi\ifmmode \check{s}\else
  \v{s}\fi{}i\ifmmode~\acute{c}\else \'{c}\fi{}}, \citenamefont {Allshouse},
  \citenamefont {Atis}, \citenamefont {Thiffeault},\ and\ \citenamefont
  {Peacock}}]{PhysRevFluids.5.054504}%
  \BibitemOpen
  \bibfield  {author} {\bibinfo {author} {\bibfnamefont {M.}~\bibnamefont
  {Filippi}}, \bibinfo {author} {\bibfnamefont {M.}~\bibnamefont {Budi\ifmmode
  \check{s}\else \v{s}\fi{}i\ifmmode~\acute{c}\else \'{c}\fi{}}}, \bibinfo
  {author} {\bibfnamefont {M.~R.}\ \bibnamefont {Allshouse}}, \bibinfo {author}
  {\bibfnamefont {S.}~\bibnamefont {Atis}}, \bibinfo {author} {\bibfnamefont
  {J.-L.}\ \bibnamefont {Thiffeault}}, \ and\ \bibinfo {author} {\bibfnamefont
  {T.}~\bibnamefont {Peacock}},\ }\href {\doibase
  10.1103/PhysRevFluids.5.054504} {\bibfield  {journal} {\bibinfo  {journal}
  {Phys. Rev. Fluids}\ }\textbf {\bibinfo {volume} {5}},\ \bibinfo {pages}
  {054504} (\bibinfo {year} {2020})}\BibitemShut {NoStop}%
\bibitem [{\citenamefont {Haller}(2015)}]{LCShaller}%
  \BibitemOpen
  \bibfield  {author} {\bibinfo {author} {\bibfnamefont {G.}~\bibnamefont
  {Haller}},\ }\href {\doibase 10.1146/annurev-fluid-010313-141322} {\bibfield
  {journal} {\bibinfo  {journal} {Annual Review of Fluid Mechanics}\ }\textbf
  {\bibinfo {volume} {47}},\ \bibinfo {pages} {137} (\bibinfo {year}
  {2015})}\BibitemShut {NoStop}%
\bibitem [{\citenamefont {Boyland}\ \emph {et~al.}(2003)\citenamefont
  {Boyland}, \citenamefont {Stremler},\ and\ \citenamefont
  {Aref}}]{BOYLAND200369}%
  \BibitemOpen
  \bibfield  {author} {\bibinfo {author} {\bibfnamefont {P.}~\bibnamefont
  {Boyland}}, \bibinfo {author} {\bibfnamefont {M.}~\bibnamefont {Stremler}}, \
  and\ \bibinfo {author} {\bibfnamefont {H.}~\bibnamefont {Aref}},\ }\href
  {\doibase https://doi.org/10.1016/S0167-2789(02)00692-9} {\bibfield
  {journal} {\bibinfo  {journal} {Physica D: Nonlinear Phenomena}\ }\textbf
  {\bibinfo {volume} {175}},\ \bibinfo {pages} {69} (\bibinfo {year}
  {2003})}\BibitemShut {NoStop}%
\bibitem [{\citenamefont {Smith}(2015)}]{smith2015point}%
  \BibitemOpen
  \bibfield  {author} {\bibinfo {author} {\bibfnamefont {S.~A.}\ \bibnamefont
  {Smith}},\ }\href@noop {} {\enquote {\bibinfo {title} {Point vortices:
  Finding periodic orbits and their topological classification},}\ } (\bibinfo
  {year} {2015}),\ \Eprint {http://arxiv.org/abs/1510.06756} {arXiv:1510.06756
  [nlin.CD]} \BibitemShut {NoStop}%
\bibitem [{\citenamefont {Stremler}\ and\ \citenamefont
  {Chen}(2007)}]{doi:10.1063/1.2772881}%
  \BibitemOpen
  \bibfield  {author} {\bibinfo {author} {\bibfnamefont {M.~A.}\ \bibnamefont
  {Stremler}}\ and\ \bibinfo {author} {\bibfnamefont {J.}~\bibnamefont
  {Chen}},\ }\href {\doibase 10.1063/1.2772881} {\bibfield  {journal} {\bibinfo
   {journal} {Physics of Fluids}\ }\textbf {\bibinfo {volume} {19}},\ \bibinfo
  {pages} {103602} (\bibinfo {year} {2007})},\ \Eprint
  {http://arxiv.org/abs/https://doi.org/10.1063/1.2772881}
  {https://doi.org/10.1063/1.2772881} \BibitemShut {NoStop}%
\bibitem [{\citenamefont {Chen}\ and\ \citenamefont
  {Stremler}(2009)}]{doi:10.1063/1.3076247}%
  \BibitemOpen
  \bibfield  {author} {\bibinfo {author} {\bibfnamefont {J.}~\bibnamefont
  {Chen}}\ and\ \bibinfo {author} {\bibfnamefont {M.~A.}\ \bibnamefont
  {Stremler}},\ }\href {\doibase 10.1063/1.3076247} {\bibfield  {journal}
  {\bibinfo  {journal} {Physics of Fluids}\ }\textbf {\bibinfo {volume} {21}},\
  \bibinfo {pages} {021701} (\bibinfo {year} {2009})},\ \Eprint
  {http://arxiv.org/abs/https://doi.org/10.1063/1.3076247}
  {https://doi.org/10.1063/1.3076247} \BibitemShut {NoStop}%
\bibitem [{\citenamefont {Stremler}\ \emph {et~al.}(2011)\citenamefont
  {Stremler}, \citenamefont {Ross}, \citenamefont {Grover},\ and\ \citenamefont
  {Kumar}}]{PhysRevLett.106.114101}%
  \BibitemOpen
  \bibfield  {author} {\bibinfo {author} {\bibfnamefont {M.~A.}\ \bibnamefont
  {Stremler}}, \bibinfo {author} {\bibfnamefont {S.~D.}\ \bibnamefont {Ross}},
  \bibinfo {author} {\bibfnamefont {P.}~\bibnamefont {Grover}}, \ and\ \bibinfo
  {author} {\bibfnamefont {P.}~\bibnamefont {Kumar}},\ }\href {\doibase
  10.1103/PhysRevLett.106.114101} {\bibfield  {journal} {\bibinfo  {journal}
  {Phys. Rev. Lett.}\ }\textbf {\bibinfo {volume} {106}},\ \bibinfo {pages}
  {114101} (\bibinfo {year} {2011})}\BibitemShut {NoStop}%
\bibitem [{\citenamefont {Grover}\ \emph {et~al.}(2012)\citenamefont {Grover},
  \citenamefont {Ross}, \citenamefont {Stremler},\ and\ \citenamefont
  {Kumar}}]{doi:10.1063/1.4768666}%
  \BibitemOpen
  \bibfield  {author} {\bibinfo {author} {\bibfnamefont {P.}~\bibnamefont
  {Grover}}, \bibinfo {author} {\bibfnamefont {S.~D.}\ \bibnamefont {Ross}},
  \bibinfo {author} {\bibfnamefont {M.~A.}\ \bibnamefont {Stremler}}, \ and\
  \bibinfo {author} {\bibfnamefont {P.}~\bibnamefont {Kumar}},\ }\href
  {\doibase 10.1063/1.4768666} {\bibfield  {journal} {\bibinfo  {journal}
  {Chaos: An Interdisciplinary Journal of Nonlinear Science}\ }\textbf
  {\bibinfo {volume} {22}},\ \bibinfo {pages} {043135} (\bibinfo {year}
  {2012})},\ \Eprint {http://arxiv.org/abs/https://doi.org/10.1063/1.4768666}
  {https://doi.org/10.1063/1.4768666} \BibitemShut {NoStop}%
\bibitem [{\citenamefont {Finn}\ and\ \citenamefont
  {Thiffeault}(2011)}]{MR2861264}%
  \BibitemOpen
  \bibfield  {author} {\bibinfo {author} {\bibfnamefont {M.~D.}\ \bibnamefont
  {Finn}}\ and\ \bibinfo {author} {\bibfnamefont {J.-L.}\ \bibnamefont
  {Thiffeault}},\ }\href {\doibase 10.1137/100791828} {\bibfield  {journal}
  {\bibinfo  {journal} {SIAM Rev.}\ }\textbf {\bibinfo {volume} {53}},\
  \bibinfo {pages} {723} (\bibinfo {year} {2011})}\BibitemShut {NoStop}%
\bibitem [{\citenamefont {Smith}\ and\ \citenamefont
  {Warrier}(2016)}]{ReynoldMixSmith}%
  \BibitemOpen
  \bibfield  {author} {\bibinfo {author} {\bibfnamefont {S.~A.}\ \bibnamefont
  {Smith}}\ and\ \bibinfo {author} {\bibfnamefont {S.}~\bibnamefont
  {Warrier}},\ }\href {\doibase 10.1063/1.4943170} {\bibfield  {journal}
  {\bibinfo  {journal} {Chaos: An Interdisciplinary Journal of Nonlinear
  Science}\ }\textbf {\bibinfo {volume} {26}},\ \bibinfo {pages} {033106}
  (\bibinfo {year} {2016})}\BibitemShut {NoStop}%
\bibitem [{\citenamefont {Thiffeault}(2018)}]{TaffyThiffeault}%
  \BibitemOpen
  \bibfield  {author} {\bibinfo {author} {\bibfnamefont {J.-L.}\ \bibnamefont
  {Thiffeault}},\ }\href {\doibase 10.1007/s00283-018-9788-4} {\bibfield
  {journal} {\bibinfo  {journal} {The Mathematical Intelligencer}\ }\textbf
  {\bibinfo {volume} {40}},\ \bibinfo {pages} {26} (\bibinfo {year}
  {2018})}\BibitemShut {NoStop}%
\bibitem [{\citenamefont {Tan}\ \emph {et~al.}(2019)\citenamefont {Tan},
  \citenamefont {Roberts}, \citenamefont {Smith}, \citenamefont {Olvera},
  \citenamefont {Arteaga}, \citenamefont {Fortini}, \citenamefont {Mitchell},\
  and\ \citenamefont {Hirst}}]{ActiveNemTopChaos2019}%
  \BibitemOpen
  \bibfield  {author} {\bibinfo {author} {\bibfnamefont {A.~J.}\ \bibnamefont
  {Tan}}, \bibinfo {author} {\bibfnamefont {E.}~\bibnamefont {Roberts}},
  \bibinfo {author} {\bibfnamefont {S.~A.}\ \bibnamefont {Smith}}, \bibinfo
  {author} {\bibfnamefont {U.~A.}\ \bibnamefont {Olvera}}, \bibinfo {author}
  {\bibfnamefont {J.}~\bibnamefont {Arteaga}}, \bibinfo {author} {\bibfnamefont
  {S.}~\bibnamefont {Fortini}}, \bibinfo {author} {\bibfnamefont {K.~A.}\
  \bibnamefont {Mitchell}}, \ and\ \bibinfo {author} {\bibfnamefont {L.~S.}\
  \bibnamefont {Hirst}},\ }\href {\doibase 10.1038/s41567-019-0600-y}
  {\bibfield  {journal} {\bibinfo  {journal} {Nature Physics}\ }\textbf
  {\bibinfo {volume} {15}},\ \bibinfo {pages} {1033} (\bibinfo {year}
  {2019})}\BibitemShut {NoStop}%
\bibitem [{\citenamefont {Budi\v{s}i\'{c}}\ and\ \citenamefont
  {Thiffeault}(2015)}]{MR3456018}%
  \BibitemOpen
  \bibfield  {author} {\bibinfo {author} {\bibfnamefont {M.}~\bibnamefont
  {Budi\v{s}i\'{c}}}\ and\ \bibinfo {author} {\bibfnamefont {J.-L.}\
  \bibnamefont {Thiffeault}},\ }\href {\doibase 10.1063/1.4927438} {\bibfield
  {journal} {\bibinfo  {journal} {Chaos}\ }\textbf {\bibinfo {volume} {25}},\
  \bibinfo {pages} {087407, 12} (\bibinfo {year} {2015})}\BibitemShut {NoStop}%
\bibitem [{\citenamefont {Nielsen}(1944)}]{MR15791}%
  \BibitemOpen
  \bibfield  {author} {\bibinfo {author} {\bibfnamefont {J.}~\bibnamefont
  {Nielsen}},\ }\href@noop {} {\bibfield  {journal} {\bibinfo  {journal}
  {Danske Vid. Selsk. Mat.-Fys. Medd.}\ }\textbf {\bibinfo {volume} {21}},\
  \bibinfo {pages} {89} (\bibinfo {year} {1944})}\BibitemShut {NoStop}%
\bibitem [{\citenamefont {Thurston}(1988)}]{MR956596}%
  \BibitemOpen
  \bibfield  {author} {\bibinfo {author} {\bibfnamefont {W.~P.}\ \bibnamefont
  {Thurston}},\ }\href {\doibase 10.1090/S0273-0979-1988-15685-6} {\bibfield
  {journal} {\bibinfo  {journal} {Bull. Amer. Math. Soc. (N.S.)}\ }\textbf
  {\bibinfo {volume} {19}},\ \bibinfo {pages} {417} (\bibinfo {year}
  {1988})}\BibitemShut {NoStop}%
\bibitem [{\citenamefont {Fathi}\ \emph {et~al.}(2012)\citenamefont {Fathi},
  \citenamefont {Laudenbach},\ and\ \citenamefont {Po\'{e}naru}}]{MR3053012}%
  \BibitemOpen
  \bibfield  {author} {\bibinfo {author} {\bibfnamefont {A.}~\bibnamefont
  {Fathi}}, \bibinfo {author} {\bibfnamefont {F.}~\bibnamefont {Laudenbach}}, \
  and\ \bibinfo {author} {\bibfnamefont {V.}~\bibnamefont {Po\'{e}naru}},\
  }\href@noop {} {\emph {\bibinfo {title} {Thurston's work on surfaces}}},\
  \bibinfo {series} {Mathematical Notes}, Vol.~\bibinfo {volume} {48}\
  (\bibinfo  {publisher} {Princeton University Press, Princeton, NJ},\ \bibinfo
  {year} {2012})\ pp.\ \bibinfo {pages} {xvi+254},\ \bibinfo {note} {translated
  from the 1979 French original by Djun M. Kim and Dan Margalit}\BibitemShut
  {NoStop}%
\bibitem [{\citenamefont {Casson}\ and\ \citenamefont
  {Bleiler}(1988)}]{MR964685}%
  \BibitemOpen
  \bibfield  {author} {\bibinfo {author} {\bibfnamefont {A.~J.}\ \bibnamefont
  {Casson}}\ and\ \bibinfo {author} {\bibfnamefont {S.~A.}\ \bibnamefont
  {Bleiler}},\ }\href {\doibase 10.1017/CBO9780511623912} {\emph {\bibinfo
  {title} {Automorphisms of surfaces after {N}ielsen and {T}hurston}}},\
  \bibinfo {series} {London Mathematical Society Student Texts}, Vol.~\bibinfo
  {volume} {9}\ (\bibinfo  {publisher} {Cambridge University Press,
  Cambridge},\ \bibinfo {year} {1988})\ pp.\ \bibinfo {pages}
  {iv+105}\BibitemShut {NoStop}%
\bibitem [{\citenamefont {Adler}\ \emph {et~al.}(1965)\citenamefont {Adler},
  \citenamefont {Konheim},\ and\ \citenamefont {McAndrew}}]{MR175106}%
  \BibitemOpen
  \bibfield  {author} {\bibinfo {author} {\bibfnamefont {R.~L.}\ \bibnamefont
  {Adler}}, \bibinfo {author} {\bibfnamefont {A.~G.}\ \bibnamefont {Konheim}},
  \ and\ \bibinfo {author} {\bibfnamefont {M.~H.}\ \bibnamefont {McAndrew}},\
  }\href {\doibase 10.2307/1994177} {\bibfield  {journal} {\bibinfo  {journal}
  {Trans. Amer. Math. Soc.}\ }\textbf {\bibinfo {volume} {114}},\ \bibinfo
  {pages} {309} (\bibinfo {year} {1965})}\BibitemShut {NoStop}%
\bibitem [{\citenamefont {Bowen}(1971)}]{MR274707}%
  \BibitemOpen
  \bibfield  {author} {\bibinfo {author} {\bibfnamefont {R.}~\bibnamefont
  {Bowen}},\ }\href {\doibase 10.2307/1995565} {\bibfield  {journal} {\bibinfo
  {journal} {Trans. Amer. Math. Soc.}\ }\textbf {\bibinfo {volume} {153}},\
  \bibinfo {pages} {401} (\bibinfo {year} {1971})}\BibitemShut {NoStop}%
\bibitem [{\citenamefont {Dinaburg}(1970)}]{MR0255765}%
  \BibitemOpen
  \bibfield  {author} {\bibinfo {author} {\bibfnamefont {E.~I.}\ \bibnamefont
  {Dinaburg}},\ }\href@noop {} {\bibfield  {journal} {\bibinfo  {journal}
  {Dokl. Akad. Nauk SSSR}\ }\textbf {\bibinfo {volume} {190}},\ \bibinfo
  {pages} {19} (\bibinfo {year} {1970})}\BibitemShut {NoStop}%
\bibitem [{\citenamefont {Aref}(1984)}]{aref_1984}%
  \BibitemOpen
  \bibfield  {author} {\bibinfo {author} {\bibfnamefont {H.}~\bibnamefont
  {Aref}},\ }\href {\doibase 10.1017/S0022112084001233} {\bibfield  {journal}
  {\bibinfo  {journal} {Journal of Fluid Mechanics}\ }\textbf {\bibinfo
  {volume} {143}},\ \bibinfo {pages} {1–21} (\bibinfo {year}
  {1984})}\BibitemShut {NoStop}%
\bibitem [{\citenamefont {Ottino}\ \emph {et~al.}(1989)\citenamefont {Ottino},
  \citenamefont {M}, \citenamefont {Crighton}, \citenamefont {Ablowitz},
  \citenamefont {Davis}, \citenamefont {Hinch}, \citenamefont {Iserles},
  \citenamefont {Ockendon},\ and\ \citenamefont
  {Olver}}]{ottino1989kinematics}%
  \BibitemOpen
  \bibfield  {author} {\bibinfo {author} {\bibfnamefont {J.}~\bibnamefont
  {Ottino}}, \bibinfo {author} {\bibfnamefont {O.}~\bibnamefont {M}}, \bibinfo
  {author} {\bibfnamefont {D.}~\bibnamefont {Crighton}}, \bibinfo {author}
  {\bibfnamefont {M.}~\bibnamefont {Ablowitz}}, \bibinfo {author}
  {\bibfnamefont {S.}~\bibnamefont {Davis}}, \bibinfo {author} {\bibfnamefont
  {E.}~\bibnamefont {Hinch}}, \bibinfo {author} {\bibfnamefont
  {A.}~\bibnamefont {Iserles}}, \bibinfo {author} {\bibfnamefont
  {J.}~\bibnamefont {Ockendon}}, \ and\ \bibinfo {author} {\bibfnamefont
  {P.}~\bibnamefont {Olver}},\ }\href
  {https://books.google.com/books?id=8OLVcbRoNSgC} {\emph {\bibinfo {title}
  {The Kinematics of Mixing: Stretching, Chaos, and Transport}}},\ Cambridge
  Texts in Applied Mathematics\ (\bibinfo  {publisher} {Cambridge University
  Press},\ \bibinfo {year} {1989})\BibitemShut {NoStop}%
\bibitem [{\citenamefont {Lanneau}\ and\ \citenamefont
  {Thiffeault}(2011{\natexlab{a}})}]{MR2828128}%
  \BibitemOpen
  \bibfield  {author} {\bibinfo {author} {\bibfnamefont {E.}~\bibnamefont
  {Lanneau}}\ and\ \bibinfo {author} {\bibfnamefont {J.-L.}\ \bibnamefont
  {Thiffeault}},\ }\href {\doibase 10.5802/aif.2599} {\bibfield  {journal}
  {\bibinfo  {journal} {Ann. Inst. Fourier (Grenoble)}\ }\textbf {\bibinfo
  {volume} {61}},\ \bibinfo {pages} {105} (\bibinfo {year}
  {2011}{\natexlab{a}})}\BibitemShut {NoStop}%
\bibitem [{\citenamefont {Lanneau}\ and\ \citenamefont
  {Thiffeault}(2011{\natexlab{b}})}]{MR2795241}%
  \BibitemOpen
  \bibfield  {author} {\bibinfo {author} {\bibfnamefont {E.}~\bibnamefont
  {Lanneau}}\ and\ \bibinfo {author} {\bibfnamefont {J.-L.}\ \bibnamefont
  {Thiffeault}},\ }\href {\doibase 10.1007/s10711-010-9551-2} {\bibfield
  {journal} {\bibinfo  {journal} {Geom. Dedicata}\ }\textbf {\bibinfo {volume}
  {152}},\ \bibinfo {pages} {165} (\bibinfo {year} {2011}{\natexlab{b}})},\
  \bibinfo {note} {supplementary material available online}\BibitemShut
  {NoStop}%
\bibitem [{\citenamefont {Hironaka}\ and\ \citenamefont
  {Kin}(2006)}]{10.2140/agt.2006.6.699}%
  \BibitemOpen
  \bibfield  {author} {\bibinfo {author} {\bibfnamefont {E.}~\bibnamefont
  {Hironaka}}\ and\ \bibinfo {author} {\bibfnamefont {E.}~\bibnamefont {Kin}},\
  }\href {\doibase 10.2140/agt.2006.6.699} {\bibfield  {journal} {\bibinfo
  {journal} {Algebraic \& Geometric Topology}\ }\textbf {\bibinfo {volume}
  {6}},\ \bibinfo {pages} {699 } (\bibinfo {year} {2006})}\BibitemShut
  {NoStop}%
\bibitem [{\citenamefont {Artin}(1947)}]{MR19087}%
  \BibitemOpen
  \bibfield  {author} {\bibinfo {author} {\bibfnamefont {E.}~\bibnamefont
  {Artin}},\ }\href {\doibase 10.2307/1969218} {\bibfield  {journal} {\bibinfo
  {journal} {Ann. of Math. (2)}\ }\textbf {\bibinfo {volume} {48}},\ \bibinfo
  {pages} {101} (\bibinfo {year} {1947})}\BibitemShut {NoStop}%
\bibitem [{\citenamefont {Penner}\ and\ \citenamefont
  {Harer}(1992)}]{MR1144770}%
  \BibitemOpen
  \bibfield  {author} {\bibinfo {author} {\bibfnamefont {R.~C.}\ \bibnamefont
  {Penner}}\ and\ \bibinfo {author} {\bibfnamefont {J.~L.}\ \bibnamefont
  {Harer}},\ }\href {\doibase 10.1515/9781400882458} {\emph {\bibinfo {title}
  {Combinatorics of train tracks}}},\ \bibinfo {series} {Annals of Mathematics
  Studies}, Vol.\ \bibinfo {volume} {125}\ (\bibinfo  {publisher} {Princeton
  University Press, Princeton, NJ},\ \bibinfo {year} {1992})\ pp.\ \bibinfo
  {pages} {xii+216}\BibitemShut {NoStop}%
\bibitem [{\citenamefont {Agol}(2011)}]{Agol2011IdealTO}%
  \BibitemOpen
  \bibfield  {author} {\bibinfo {author} {\bibfnamefont {I.}~\bibnamefont
  {Agol}},\ }\href {\doibase http://dx.doi.org/10.1090/conm/560} {\bibfield
  {journal} {\bibinfo  {journal} {Contemporary Mathematics}\ }\textbf {\bibinfo
  {volume} {560}} (\bibinfo {year} {2011}),\
  http://dx.doi.org/10.1090/conm/560}\BibitemShut {NoStop}%
\bibitem [{\citenamefont {Birman}(1975)}]{MR0425944}%
  \BibitemOpen
  \bibfield  {author} {\bibinfo {author} {\bibfnamefont {J.~S.}\ \bibnamefont
  {Birman}},\ }\href@noop {} {\emph {\bibinfo {title} {Erratum: ``{\it
  {B}raids, links, and mapping class groups}'' ({A}nn. of {M}ath. {S}tudies,
  {N}o. 82, {P}rinceton {U}niv. {P}ress, {P}rinceton, {N}. {J}., 1974)}}}\
  (\bibinfo  {publisher} {Princeton University Press, Princeton, N. J.;
  University of Tokyo Press, Toyko},\ \bibinfo {year} {1975})\ p.~\bibinfo
  {pages} {1},\ \bibinfo {note} {based on lecture notes by James
  Cannon}\BibitemShut {NoStop}%
\bibitem [{\citenamefont {Kassel}\ \emph {et~al.}(2008)\citenamefont {Kassel},
  \citenamefont {Dodane},\ and\ \citenamefont {Turaev}}]{kassel2008braid}%
  \BibitemOpen
  \bibfield  {author} {\bibinfo {author} {\bibfnamefont {C.}~\bibnamefont
  {Kassel}}, \bibinfo {author} {\bibfnamefont {O.}~\bibnamefont {Dodane}}, \
  and\ \bibinfo {author} {\bibfnamefont {V.}~\bibnamefont {Turaev}},\ }\href
  {https://books.google.com/books?id=y6Cox3XjdroC} {\emph {\bibinfo {title}
  {Braid Groups}}},\ Graduate Texts in Mathematics\ (\bibinfo  {publisher}
  {Springer New York},\ \bibinfo {year} {2008})\BibitemShut {NoStop}%
\bibitem [{\citenamefont {Bellingeri}(2001)}]{Bellingeri2001OnPO}%
  \BibitemOpen
  \bibfield  {author} {\bibinfo {author} {\bibfnamefont {P.}~\bibnamefont
  {Bellingeri}},\ }\href@noop {} {\bibfield  {journal} {\bibinfo  {journal}
  {Journal of Algebra}\ }\textbf {\bibinfo {volume} {274}},\ \bibinfo {pages}
  {543} (\bibinfo {year} {2001})}\BibitemShut {NoStop}%
\bibitem [{\citenamefont {Dehornoy}(2008)}]{MR2462117}%
  \BibitemOpen
  \bibfield  {author} {\bibinfo {author} {\bibfnamefont {P.}~\bibnamefont
  {Dehornoy}},\ }\href {\doibase 10.1016/j.dam.2007.12.009} {\bibfield
  {journal} {\bibinfo  {journal} {Discrete Appl. Math.}\ }\textbf {\bibinfo
  {volume} {156}},\ \bibinfo {pages} {3091} (\bibinfo {year}
  {2008})}\BibitemShut {NoStop}%
\bibitem [{\citenamefont {Sergiescu}(1993)}]{Sergiescu1993}%
  \BibitemOpen
  \bibfield  {author} {\bibinfo {author} {\bibfnamefont {V.}~\bibnamefont
  {Sergiescu}},\ }\href {http://eudml.org/doc/174582} {\bibfield  {journal}
  {\bibinfo  {journal} {Mathematische Zeitschrift}\ }\textbf {\bibinfo {volume}
  {214}},\ \bibinfo {pages} {477} (\bibinfo {year} {1993})}\BibitemShut
  {NoStop}%
\bibitem [{\citenamefont {Mohar}\ and\ \citenamefont
  {Thomassen}(2001)}]{mohar2001graphs}%
  \BibitemOpen
  \bibfield  {author} {\bibinfo {author} {\bibfnamefont {B.}~\bibnamefont
  {Mohar}}\ and\ \bibinfo {author} {\bibfnamefont {C.}~\bibnamefont
  {Thomassen}},\ }\href@noop {} {\emph {\bibinfo {title} {Graphs on
  Surfaces}}},\ Johns Hopkins Studies in Nineteenth C Architecture Series\
  (\bibinfo  {publisher} {Johns Hopkins University Press},\ \bibinfo {year}
  {2001})\BibitemShut {NoStop}%
\bibitem [{\citenamefont {{\v{S}}ir{\'a}{\v{n}}}(2006)}]{RegularMapsSurvey}%
  \BibitemOpen
  \bibfield  {author} {\bibinfo {author} {\bibfnamefont {J.}~\bibnamefont
  {{\v{S}}ir{\'a}{\v{n}}}},\ }in\ \href@noop {} {\emph {\bibinfo {booktitle}
  {Topics in Discrete Mathematics}}},\ \bibinfo {editor} {edited by\ \bibinfo
  {editor} {\bibfnamefont {M.}~\bibnamefont {Klazar}}, \bibinfo {editor}
  {\bibfnamefont {J.}~\bibnamefont {Kratochv{\'i}l}}, \bibinfo {editor}
  {\bibfnamefont {M.}~\bibnamefont {Loebl}}, \bibinfo {editor} {\bibfnamefont
  {J.}~\bibnamefont {Matou{\v{s}}ek}}, \bibinfo {editor} {\bibfnamefont
  {P.}~\bibnamefont {Valtr}}, \ and\ \bibinfo {editor} {\bibfnamefont
  {R.}~\bibnamefont {Thomas}}}\ (\bibinfo  {publisher} {Springer Berlin
  Heidelberg},\ \bibinfo {address} {Berlin, Heidelberg},\ \bibinfo {year}
  {2006})\ pp.\ \bibinfo {pages} {591--609}\BibitemShut {NoStop}%
\bibitem [{\citenamefont {Finn}\ and\ \citenamefont
  {Thiffeault}(2007)}]{MR2299637}%
  \BibitemOpen
  \bibfield  {author} {\bibinfo {author} {\bibfnamefont {M.~D.}\ \bibnamefont
  {Finn}}\ and\ \bibinfo {author} {\bibfnamefont {J.-L.}\ \bibnamefont
  {Thiffeault}},\ }\href {\doibase 10.1137/060659636} {\bibfield  {journal}
  {\bibinfo  {journal} {SIAM J. Appl. Dyn. Syst.}\ }\textbf {\bibinfo {volume}
  {6}},\ \bibinfo {pages} {79} (\bibinfo {year} {2007})}\BibitemShut {NoStop}%
\bibitem [{\citenamefont {Farb}\ and\ \citenamefont
  {Margalit}(2012)}]{MR2850125}%
  \BibitemOpen
  \bibfield  {author} {\bibinfo {author} {\bibfnamefont {B.}~\bibnamefont
  {Farb}}\ and\ \bibinfo {author} {\bibfnamefont {D.}~\bibnamefont
  {Margalit}},\ }\href@noop {} {\emph {\bibinfo {title} {A primer on mapping
  class groups}}},\ \bibinfo {series} {Princeton Mathematical Series},
  Vol.~\bibinfo {volume} {49}\ (\bibinfo  {publisher} {Princeton University
  Press, Princeton, NJ},\ \bibinfo {year} {2012})\ pp.\ \bibinfo {pages}
  {xiv+472}\BibitemShut {NoStop}%
\bibitem [{\citenamefont {Bestvina}\ and\ \citenamefont
  {Handel}(1995)}]{MR1308491}%
  \BibitemOpen
  \bibfield  {author} {\bibinfo {author} {\bibfnamefont {M.}~\bibnamefont
  {Bestvina}}\ and\ \bibinfo {author} {\bibfnamefont {M.}~\bibnamefont
  {Handel}},\ }\href {\doibase 10.1016/0040-9383(94)E0009-9} {\bibfield
  {journal} {\bibinfo  {journal} {Topology}\ }\textbf {\bibinfo {volume}
  {34}},\ \bibinfo {pages} {109} (\bibinfo {year} {1995})}\BibitemShut
  {NoStop}%
\bibitem [{\citenamefont {Dynnikov}(2002)}]{MR1918864}%
  \BibitemOpen
  \bibfield  {author} {\bibinfo {author} {\bibfnamefont {I.~A.}\ \bibnamefont
  {Dynnikov}},\ }\href {\doibase 10.1070/RM2002v057n03ABEH000519} {\bibfield
  {journal} {\bibinfo  {journal} {Uspekhi Mat. Nauk}\ }\textbf {\bibinfo
  {volume} {57}},\ \bibinfo {pages} {151} (\bibinfo {year} {2002})}\BibitemShut
  {NoStop}%
\bibitem [{\citenamefont {Hall}\ and\ \citenamefont
  {Yurtta\c{s}}(2009)}]{MR2512607}%
  \BibitemOpen
  \bibfield  {author} {\bibinfo {author} {\bibfnamefont {T.}~\bibnamefont
  {Hall}}\ and\ \bibinfo {author} {\bibfnamefont {S.~O.}\ \bibnamefont
  {Yurtta\c{s}}},\ }\href {\doibase 10.1016/j.topol.2009.01.005} {\bibfield
  {journal} {\bibinfo  {journal} {Topology Appl.}\ }\textbf {\bibinfo {volume}
  {156}},\ \bibinfo {pages} {1554} (\bibinfo {year} {2009})}\BibitemShut
  {NoStop}%
\bibitem [{\citenamefont {Moussafir}(2006)}]{MR2381961}%
  \BibitemOpen
  \bibfield  {author} {\bibinfo {author} {\bibfnamefont {J.-O.}\ \bibnamefont
  {Moussafir}},\ }\href {\doibase 10.1007/s11853-007-0004-x} {\bibfield
  {journal} {\bibinfo  {journal} {Funct. Anal. Other Math.}\ }\textbf {\bibinfo
  {volume} {1}},\ \bibinfo {pages} {37} (\bibinfo {year} {2006})}\BibitemShut
  {NoStop}%
\bibitem [{\citenamefont {Roberts}\ \emph {et~al.}(2019)\citenamefont
  {Roberts}, \citenamefont {Sindi}, \citenamefont {Smith},\ and\ \citenamefont
  {Mitchell}}]{E-tec}%
  \BibitemOpen
  \bibfield  {author} {\bibinfo {author} {\bibfnamefont {E.}~\bibnamefont
  {Roberts}}, \bibinfo {author} {\bibfnamefont {S.}~\bibnamefont {Sindi}},
  \bibinfo {author} {\bibfnamefont {S.~A.}\ \bibnamefont {Smith}}, \ and\
  \bibinfo {author} {\bibfnamefont {K.~A.}\ \bibnamefont {Mitchell}},\ }\href
  {\doibase 10.1063/1.5045060} {\bibfield  {journal} {\bibinfo  {journal}
  {Chaos: An Interdisciplinary Journal of Nonlinear Science}\ }\textbf
  {\bibinfo {volume} {29}},\ \bibinfo {pages} {013124} (\bibinfo {year}
  {2019})}\BibitemShut {NoStop}%
\bibitem [{\citenamefont {Butkovi\v{c}}(2010)}]{MR2681232}%
  \BibitemOpen
  \bibfield  {author} {\bibinfo {author} {\bibfnamefont {P.}~\bibnamefont
  {Butkovi\v{c}}},\ }\href {\doibase 10.1007/978-1-84996-299-5} {\emph
  {\bibinfo {title} {Max-linear systems: theory and algorithms}}},\ Springer
  Monographs in Mathematics\ (\bibinfo  {publisher} {Springer-Verlag London,
  Ltd., London},\ \bibinfo {year} {2010})\ pp.\ \bibinfo {pages}
  {xviii+272}\BibitemShut {NoStop}%
\bibitem [{\citenamefont {Heidergott}\ \emph {et~al.}(2006)\citenamefont
  {Heidergott}, \citenamefont {Oldser},\ and\ \citenamefont {van~der
  Woude}}]{MR2188299}%
  \BibitemOpen
  \bibfield  {author} {\bibinfo {author} {\bibfnamefont {B.}~\bibnamefont
  {Heidergott}}, \bibinfo {author} {\bibfnamefont {G.~J.}\ \bibnamefont
  {Oldser}}, \ and\ \bibinfo {author} {\bibfnamefont {J.}~\bibnamefont {van~der
  Woude}},\ }\href@noop {} {\emph {\bibinfo {title} {Max plus at work}}},\
  Princeton Series in Applied Mathematics\ (\bibinfo  {publisher} {Princeton
  University Press, Princeton, NJ},\ \bibinfo {year} {2006})\ pp.\ \bibinfo
  {pages} {xii+213},\ \bibinfo {note} {modeling and analysis of synchronized
  systems: a course on max-plus algebra and its applications}\BibitemShut
  {NoStop}%
\bibitem [{\citenamefont {Fomin}\ \emph {et~al.}(2008)\citenamefont {Fomin},
  \citenamefont {Shapiro},\ and\ \citenamefont
  {Thurston}}]{ClusterAlgebraDThurstonI}%
  \BibitemOpen
  \bibfield  {author} {\bibinfo {author} {\bibfnamefont {S.}~\bibnamefont
  {Fomin}}, \bibinfo {author} {\bibfnamefont {M.}~\bibnamefont {Shapiro}}, \
  and\ \bibinfo {author} {\bibfnamefont {D.}~\bibnamefont {Thurston}},\ }\href
  {\doibase 10.1007/s11511-008-0030-7} {\bibfield  {journal} {\bibinfo
  {journal} {Acta Mathematica}\ }\textbf {\bibinfo {volume} {201}},\ \bibinfo
  {pages} {83} (\bibinfo {year} {2008})}\BibitemShut {NoStop}%
\bibitem [{\citenamefont {Fomin}\ and\ \citenamefont
  {Thurston}(2018)}]{fomin2018cluster}%
  \BibitemOpen
  \bibfield  {author} {\bibinfo {author} {\bibfnamefont {S.}~\bibnamefont
  {Fomin}}\ and\ \bibinfo {author} {\bibfnamefont {D.}~\bibnamefont
  {Thurston}},\ }\href@noop {} {\emph {\bibinfo {title} {Cluster Algebras and
  Triangulated Surfaces: Lambda lengths}}},\ \bibinfo {series} {Cluster
  Algebras and Triangulated Surfaces: Lambda Lengths}\ No.\ \bibinfo {number}
  {pt. 2}\ (\bibinfo  {publisher} {American Mathematical Society},\ \bibinfo
  {year} {2018})\BibitemShut {NoStop}%
\bibitem [{\citenamefont {Sanchez}\ \emph {et~al.}(2012)\citenamefont
  {Sanchez}, \citenamefont {Chen}, \citenamefont {DeCamp}, \citenamefont
  {Heymann},\ and\ \citenamefont {Dogic}}]{DogicANMT}%
  \BibitemOpen
  \bibfield  {author} {\bibinfo {author} {\bibfnamefont {T.}~\bibnamefont
  {Sanchez}}, \bibinfo {author} {\bibfnamefont {D.~T.~N.}\ \bibnamefont
  {Chen}}, \bibinfo {author} {\bibfnamefont {S.~J.}\ \bibnamefont {DeCamp}},
  \bibinfo {author} {\bibfnamefont {M.}~\bibnamefont {Heymann}}, \ and\
  \bibinfo {author} {\bibfnamefont {Z.}~\bibnamefont {Dogic}},\ }\href
  {\doibase 10.1038/nature11591} {\bibfield  {journal} {\bibinfo  {journal}
  {Nature}\ }\textbf {\bibinfo {volume} {491}},\ \bibinfo {pages} {431}
  (\bibinfo {year} {2012})}\BibitemShut {NoStop}%
\bibitem [{\citenamefont {Shendruk}\ \emph {et~al.}(2017)\citenamefont
  {Shendruk}, \citenamefont {Doostmohammadi}, \citenamefont {Thijssen},\ and\
  \citenamefont {Yeomans}}]{C6SM02310J}%
  \BibitemOpen
  \bibfield  {author} {\bibinfo {author} {\bibfnamefont {T.~N.}\ \bibnamefont
  {Shendruk}}, \bibinfo {author} {\bibfnamefont {A.}~\bibnamefont
  {Doostmohammadi}}, \bibinfo {author} {\bibfnamefont {K.}~\bibnamefont
  {Thijssen}}, \ and\ \bibinfo {author} {\bibfnamefont {J.~M.}\ \bibnamefont
  {Yeomans}},\ }\href {\doibase 10.1039/C6SM02310J} {\bibfield  {journal}
  {\bibinfo  {journal} {Soft Matter}\ }\textbf {\bibinfo {volume} {13}},\
  \bibinfo {pages} {3853} (\bibinfo {year} {2017})}\BibitemShut {NoStop}%
\bibitem [{\citenamefont {Impellizieri}(2016)}]{impellizieri2016domino}%
  \BibitemOpen
  \bibfield  {author} {\bibinfo {author} {\bibfnamefont {F.}~\bibnamefont
  {Impellizieri}},\ }\href@noop {} {\enquote {\bibinfo {title} {Domino tilings
  of the torus},}\ } (\bibinfo {year} {2016}),\ \Eprint
  {http://arxiv.org/abs/1601.06354} {arXiv:1601.06354 [math.CO]} \BibitemShut
  {NoStop}%
\end{thebibliography}%

\section{Appendices}

\subsection{Other Lattice Graphs}
\label{Appendix:OtherLatticeGraphs}

In section~\ref{Sec:MixEff} we introduced a method for calculating the topological entropy and TEPO for braids defined on any choice of orientable surface and embedded graph.  In section~\ref{Sec:Example}, we highlighted how this works for the 2 point regular embedding of the square lattice on the torus.  For the max TEPO search results in section~\ref{Sec:Search}, we included five more torus models for lattice braids.  Here we provide some of the information needed to work through these other five examples.

First of all, the code used to execute the max TEPO search for each of the six cases can be found at GITHUB (DOI: \url{https://doi.org/10.5281/zenodo.4779103}).  There is a separate Jupyter notebook, written in python, for each case.  These include a careful explanation of how to apply our procedure.

	\begin{figure}[htbp]
		\center
		\includegraphics[width = \linewidth]{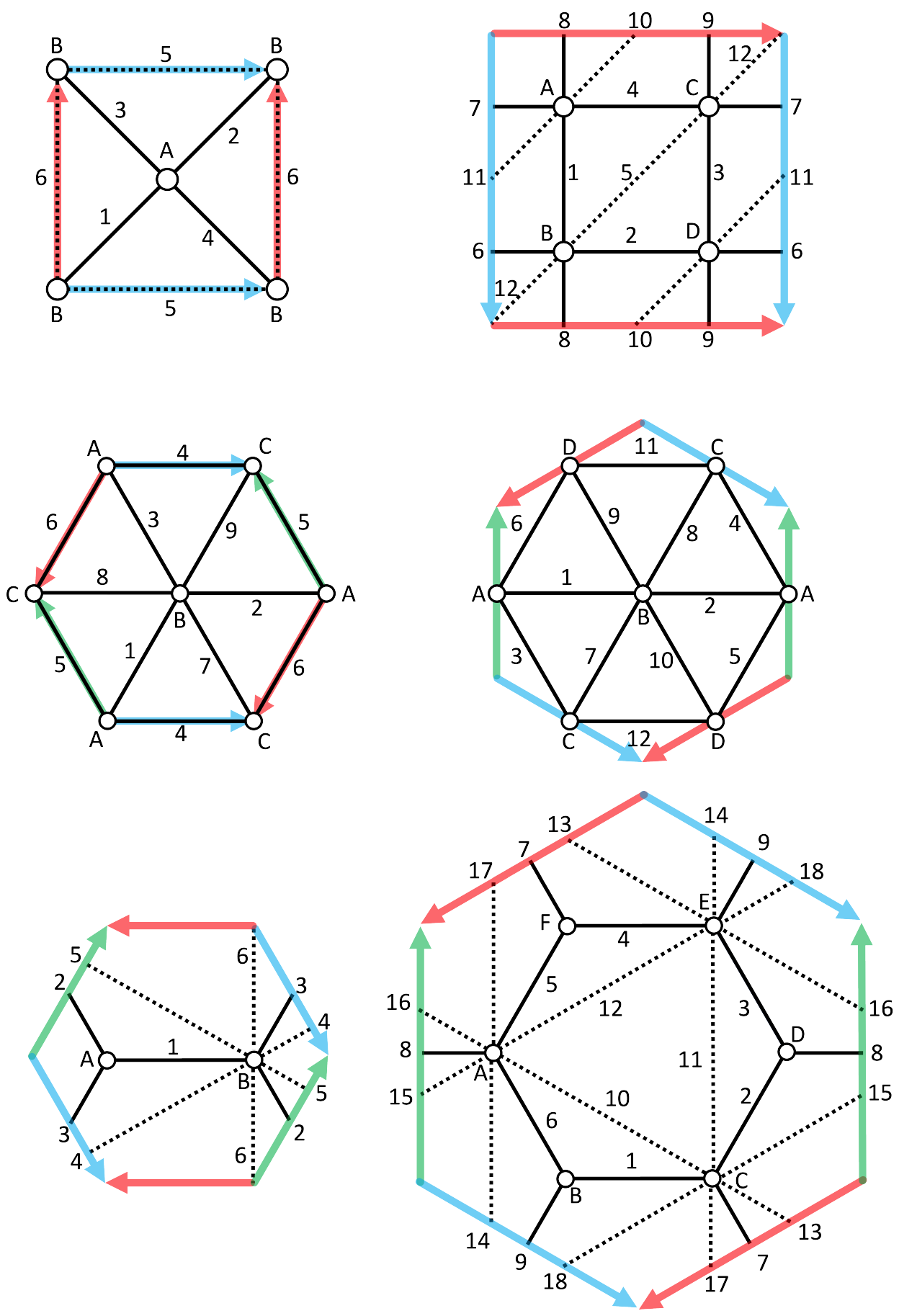}
		\caption{Torus models with labeled graph edges (solid black lines) and labeled extra edges (dashed black lines) to form triangulations.  Top: square lattice, 2 point (left) and 4 point (right), Middle: triangular lattice, 3 point (left) and 4 point (right), and Bottom: 2 point (left) and 6 point (right).}
		\label{Fig:AllCasesTriangulations}
	\end{figure}
	
\begin{table*}[t]
  \centering
  \begin{tabular}{|c|c||c|c|c|c|c|c|}
	 \hline
	 \multicolumn{2}{|c||}{Lattice Types} & \multicolumn{2}{c|}{Square} & \multicolumn{2}{c|}{Triangular} & \multicolumn{2}{c|}{Hexagonal} \\
	 \hline
	 \multicolumn{2}{|c||}{Torus Model} & 2 pt. & 4 pt. & 3 pt. & 4 pt. & 2 pt. & 6 pt. \\
	 
	  \hline
	  \hline
	 \multicolumn{2}{|c||}{\makecell{Example \\ Generator}} & $\sigma_2$ & $\sigma_1$ & $\sigma_2$ & $\sigma_2$ & $\sigma_1$ & $\sigma_1$ \\
	 \hline

	 \multicolumn{2}{|c||}{\makecell{Flip \\ Sequence}} &
	 \makecell{$6,3,5,6',$ \\ $1,5',6'',4,5''$} & \makecell{$4,6,8,10,$ \\ $12,7,2,8'$} & \makecell{$1,3,6,9,$ \\ $4,8,1',3'$} & \makecell{$1,5,8,3,$ \\ $9,6,7,1'$} & \makecell{$2,4,3,6,$ \\ $5,2',3',4'$} & \makecell{$11,17,18,6,7,9,$ \\ $13,15,2,9',11',17'$} \\
	 \hline

  \end{tabular}
  \caption{Flip Sequences}
  \label{tab:FlipSequences}
\end{table*}	

Instead of replicating the full analysis for each of the five extra cases, we will relay enough information to reconstruct the flip sequences, and therefore the intersection coordinate update rules, for a single generator.  The update rules for each of the other generators and their inverses can be found after conjugating by the appropriate graph symmetries (rotations, translations, and reflections).  The action of these symmetries on the intersection coordinates are simple enough to find through inspection.  The only difficulty arises for graphs that are not already a triangulation, and the graph symmetry doesn't automatically extend to a symmetry of the triangulation.  In these cases, the usual permutation of intersection coordinates must also be combined with triangulation flips (and the attendant coordinate updates) so as to map the whole triangulation back onto itself after the symmetry.

The minimal information needed to specify a flip sequence consists of the triangulation with labeled edges (representing the indices for the intersection coordinates), and a time-ordered sequence of edges to be flipped.  The labeled triangulations for each of the six torus models we have considered are shown in fig.~\ref{Fig:AllCasesTriangulations}.  The edges that are part of the original graph are shown as black lines, while the extra edges needed to create a triangulation are shown as dotted black lines.  The time-ordered (left to right) set of edges in the flip sequence are shown in table~\ref{tab:FlipSequences}.  Each edge in this sequence represents the edge that will be flipped in a Whitehead move.  The labeling convention is for the new edge to take the index of the flipped edge and accrue a prime (e.g. $4  \rightarrow 4'$).  The flip sequence is to be followed by a geometric movement of the two points (bounding the edge given by the index of the example generator listed in the table), which realizes the topological swap.

\subsection{General Square Lattice}
\label{Appendix:gensq}

Performing this analysis - recording graph symmetries, constructing the flip chart, and encoding the intersection coordinate update rules for each generator - takes time, and it is natural to see how much can be automated.  Ideally, given a surface graph, we could apply the ideas referenced in fig.~\ref{Fig:FlipPlanning} to each edge in the graph, even when the graph has no symmetries.  The only difficulty arises from graphs that are small enough to include pairs of points connected by more than one edge.  Indeed this difficulty is present in five of our six examples (see fig.~\ref{Fig:AllCasesTriangulations}).  For larger lattice graphs, this algorithm works well.  We have implemented this idea for square lattice graphs with $M$ rows and $N$ columns, $M,N \geq 3$.  This code, with example usage, can also be found at GITHUB (DOI: \url{https://doi.org/10.5281/zenodo.4779103}).

Here, we would like to highlight the combinatorial complexity of exhaustively searching for max TEPO braids defined on $M\times M$ square lattice graphs, as M increases.  As a reminder, braid operations - what we build our braid words out of - are sets of braid generators that pairwise commute, and therefore can be executed simultaneously. For simplicity, we consider maximal braid operations where, with $M$ even, every vertex is paired up in a generator. From the perspective of graph theory, where lattice points are vertices and graph edges represent possible generator switches, a maximal braid operation is a complete matching on the graph, along with a CW/CCW designation for each edge in the matching.  We will refer to a complete matching as a braid operation template, or template for short.

To get a feel for the combinatorial complexity of braid operations for square lattices, there are $272$ unique templates for the $4\times 4$ square lattice on a torus. Again, these templates only set which vertices each generator acts on; we must also consider whether the generator is acting in a clockwise or counterclockwise direction. For this $4\times 4$ lattice, a template is comprised of $8$ generators, so we now have $272(2^8)=69632$ possible maximal operations. 

To simplify things, we can separate these templates into subgroups based on the flux of a matching.  The flux, part of the machinery used to enumerate perfect matchings on a torus~\cite{impellizieri2016domino}, is an ordered pair of integers that describes how the graph edges cross over the fundamental domain of the torus in each direction. The first value, the horizontal flux, is the number of graph edges crossing over the left/right fundamental domain edge, counted positively if they are in an odd row (starting from the top), and counted negatively if in an even row. The second value, the vertical flux, analogously counts the number of graph edges crossing over the top/bottom fundamental domain edge; again counted positively if in an odd column (starting from the left), and negatively if in an even column. 

For the 4 by 4 square lattice on a torus, there are $132$ maximal operations in the $(0,0)$ flux group, $32$ operations in each of the four $(0,\pm 1),(\pm 1,0),$ flux groups, $2$ operations in each of the $(\pm 1,1),(\pm 1,-1)$ flux groups, and $1$ operation in each of the $(0,\pm 2),(\pm 2,0)$ flux groups. Notably, our candidate maximum TEPO braid is composed of braid operations with templates exclusively from this last set of flux groups.  Indeed, for the $2M \times 2M$ case, the analogue of our max TEPO braid is created from the $(0,\pm M)$ and $(\pm M,0)$ flux groups, which have only one representative each.  This points to the uniqueness of our candidate max TEPO braid.

For an arbitrarily sized square lattice, there are far too many potential operations for us to be able to verify the max TEPO braid through the brute force method of testing all possible combinations. For the 4 by 4 square lattice alone, there are $69632$ maximal operations, so if we only tested braids of length $4$, we would need to calculate the TEPO for $69632^4=2.35\times 10^{19}$ braids. If we increased the size of the lattice or the length of the braid, this number would grow exponentially. As noted, our candidate maximum TEPO braid has operations entirely in the $(0,\pm 2),(\pm 2,0)$ flux groups, so creating braids from specific flux groups may be a good way to approach this problem. This strategy reduces the complexity of our search, but is still combinatorially challenging. For the $4 \times 4$ square lattice, we could choose to look at braids composed entirely of operations pulled from the $(0,\pm 2),(\pm 2,0)$ flux groups. In this case, there are $4$ operation templates, and $4(2^8)=1024$ maximal operations. So, if we only check braids of length $4$, we must calculate the TEPO for $1024^4=1.1\times10^{12}$ braids. This strategy has significantly reduced the scale of our calculations, but on its own will not solve the combinatorial challenges of a larger lattice. Flux groups can help us organize braid templates and potentially reduce some of the combinatorial difficulties inherent in a exhaustive search for max TEPO braids.

\end{document}